\documentclass[12pt,reqno]{article}

\usepackage[english]{babel}		

\usepackage[latin1]{inputenc}
\usepackage{amscd}
\usepackage{amsfonts}
\usepackage{amsgen}
\usepackage{amsmath}
\usepackage{amssymb}
\usepackage{amstext}
\usepackage{amsthm}
\usepackage{graphicx}
\usepackage{latexsym}
\usepackage{epsfig}
\usepackage{wrapfig}
\usepackage[all,knot,arc,import,poly]{xy}
\usepackage[usenames]{color}
\usepackage[mathscr]{eucal}
\usepackage{indentfirst}
\usepackage[T1]{fontenc}
\usepackage{setspace}
\usepackage{lscape}
\usepackage{tikz}
\usepackage{pgfplots}
\pgfplotsset{compat=1.13}
\usepackage{braids}
\usepackage[pdftex,                %
implicit=true,
bookmarks=true,
breaklinks=true,
bookmarksnumbered = true,
colorlinks= true,
linkcolor= blue,
urlcolor= green,
anchorcolor = yellow,
citecolor=blue,
]{hyperref}%
\everymath{\displaystyle}
\usepackage{comment}

\newtheorem{Def}{Definition}[section]
\newtheorem{Teo}[Def]{Theorem}
\newtheorem{Prop}[Def]{Proposition}
\newtheorem{Lema}[Def]{Lemma}

\newtheorem{Obs}[Def]{Remark}

\newtheorem{Cor}[Def]{Corollary}

\newcommand{\che}{\left\{}
\newcommand{\chd}{\right\}}
\newcommand{\pare}{\left(}
\newcommand{\pard}{\right)}

\newcommand{\R}{\mathbb{R}}
\newcommand{\N}{\mathbb{N}}

\newcommand{\Z}{\mathbb{Z}}
\newcommand{\aut}{\operatorname{Aut}}

\newcommand{\mdc}{\operatorname{gcd}}

\newcommand{\comp}[1]{\noindent\textcolor{blue}{\textbf{!!!P!!!~#1}}}

\begin{document}

\title{Virtual braid groups, virtual twin groups and crystallographic groups}



\author{OSCAR~OCAMPO~\\
Departamento de Matem\'atica - Instituto de Matem\'atica e Estat\'istica,\\
Universidade Federal da Bahia,\\
CEP:~40170-110 - Salvador - Ba - Brazil.\\
e-mail:~\url{oscaro@ufba.br}\vspace*{4mm}\\
PAULO CESAR CERQUEIRA DOS SANTOS J\'{U}NIOR ~\\
%
Departamento de Matem\'atica - Instituto de Matem\'atica e Estat\'istica,\\
Universidade Federal da Bahia,\\
CEP:~40170-110 - Salvador - Ba - Brazil.\\
e-mail:~\url{pcesarmath@gmail.com}
}

\maketitle



\begin{abstract}

Let $n\ge 2$. Let $VB_n$ (resp. $VP_n$) be the virtual braid group (resp. the pure virtual braid group), and let $VT_n$ (resp. $PVT_n$) be the virtual twin group (resp. the pure virtual twin group). 
Let $\Pi$ be one of the following quotients: $VB_n/\Gamma_2(VP_n)$ or $VT_n/\Gamma_2(PVT_n)$ where $\Gamma_2(H)$ is the commutator subgroup of $H$. 
In this paper, we show that $\Pi$ is a crystallographic group and we characterize the elements of finite order and the conjugacy classes of elements in $\Pi$. 
Furthermore, we realize explicitly some Bieberbach groups and infinite virtually cyclic groups in $\Pi$. 
Finally, we also study other braid-like groups (welded, unrestricted, flat virtual, flat welded and Gauss virtual braid group) modulo the respective commutator subgroup in each case.


\end{abstract}

\let\thefootnote\relax

\section{Introduction}

There are several generalizations of the Artin braid group $B_n$, both from a geometric and algebraic point of view, e.g. the mapping class groups, surface braid groups, virtual braid groups, Artin-Tits groups, Garside groups, twin and virtual twin groups, among others, and constitute an active line of research. 
In this paper, we are interested in studying some quotients of the virtual braid group $VB_n$ and the virtual twin group $VT_n$. 

 Let $S_n$ be the symmetric group defined over a set $X$ of $n$ symbols. Given $\sigma,\tau\in S_n$ we shall use in this paper the operation given by $\sigma\tau(x)=\tau(\sigma(x))$ for all $x$ in the set $X$.

The virtual braid group is the natural companion to the category of virtual knots, just as the Artin braid group is to usual knots and links. 
We note that a virtual knot diagram is like a classical knot diagram with one extra type of crossing, called a virtual crossing. The virtual braid groups have interpretations in terms of diagrams, see~\cite{K2007},~\cite{K1999} and~\cite{V2001}.
The notion of virtual knots and links was introduced by Kauffman together with virtual braids in \cite{K1999}, and since then it has drawn the attention of several researchers. 
For instance, for the case of virtual braids, its homological properties were studied in~\cite{V2001}; properties and structural aspects were showed in \cite{B2004}, \cite{BB2009}; and in~\cite{BP2020} it was determined all possible homomorphisms from $VB_n$ to $S_m$, from $S_n$ to $VB_m$ and from $VB_n$ to $VB_m$ for $n,m\in\N$, $n\ge 5$, $m\ge 2$ and $n\ge m$.


The virtual twin group $VT_n$ was introduced in~\cite{BSV2019} as an abstract generalization of twin groups, together with the pure virtual twin group, $PVT_n$ which is defined as the kernel of the natural surjection from $VT_n$ onto $S_n$. 
Then, in \cite{NNS2020}, the virtual twin group $VT_n$ and the pure virtual twin group $PVT_n$ were studied with more details. 
The virtual twin group $VT_n$ contains the twin group $T_n$ and the symmetric group $S_n$ and, similarly to virtual braid groups, they have a geometrical interpretation. 
In \cite{37} the authors proved Alexander's and Markov's theorems for virtual doodles.
The virtual braid and virtual twin groups (also some of their quotients) have various applications in algebra and topology, and they have a geometric interpretation that is often very useful when studying applications of them on their structural properties. 

Many quotients of virtual braid groups are also very interesting, in particular the welded braid group $WB_n$, the unrestricted virtual braid group $UVB_n$, the flat virtual braid group $FVB_n$, the flat welded braid groups $FWB_n$, and the Gauss virtual braid group $GVB_n$. For more details about these groups see \cite{BBD2015}, \cite{D2017}, \cite{K2007} and \cite{K1999}. 
As doing for virtual braid groups, it is very natural to define similar quotients for virtual twin groups. 
For instance, the welded twin groups were defined in \cite{BSV2019}. 
In the diagram below, we also define the \textit{Gauss virtual twin groups}.

\begin{Obs}
Let $G$ be a group and let $N$ be a normal subgroup of $G$. 
By an abuse of notation, in this paper we use the same symbols for elements in $G$ and its equivalence classes in the quotient $G/N$. 

Let $\Gamma_2(N)$ denote the commutator subgroup of $N$. 
In particular, the abuse of notation mentioned above will be used several times in the case of $G$ and the quotient $G/\Gamma_2(N)$. 
We note that it also applies for the abelianization of $N$.

\end{Obs}

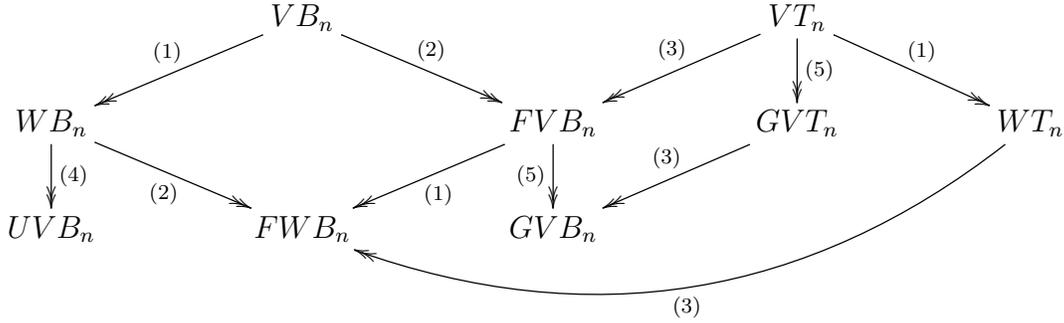
\begin{figure}[h] 
 \[ 
    \xymatrix@C=4.5pc{
  & VB_n 
  \ar@{->>}[dl]_{(1)} 
  \ar@{->>}[dr]^-{(2)}
    &    & 
  VT_n \ar@{->>}[d]^-{(5)}
  \ar@{->>}[dl]_-{(3)}
  \ar@{->>}[dr]^-{(1)} & \\
WB_n \ar@{->>}[d]^-{(4)}
\ar@{->>}[dr]_-{(2)} &  & 
FVB_n 
\ar@{->>}[d]_-{(5)}
\ar@{->>}[dl]^-{(1)} & 
GVT_n 
\ar@{->>}[dl]_-{(3)}  & 
WT_n \ar@/^1.5cm/@{->>}[dlll]^{(3)}
 \\
UVB_n &  FWB_n  & GVB_n &  &   
}
 \]
\caption{Quotients of virtual braid and virtual twin groups}\label{quocientes}
\end{figure}

Presentations for $VB_n$ and $VT_n$ may be found in Theorems~\ref{apvbn} and~\ref{apvtn}, respectively. 
In Figure~\ref{quocientes} we show in a diagram some quotients of virtual braid and virtual twin groups. 
The surjective homomorphism (directed arrow) between two groups means that the target is obtained, from the other group, just adding to the presentation of the former one the set of relations indicated by the number given in the arrow, as follows:
\begin{itemize}
    \item[(1)] $\rho_i\sigma_{i+1}\sigma_i = \sigma_{i+1}\sigma_i\rho_{i+1}$, for $i=1,\ldots, n-2$,
\item[(2)] $\sigma_i^2=1$, for $i=1,\ldots, n-1$,
\item[(3)] $\sigma_i\sigma_{i+1}\sigma_i = \sigma_{i+1}\sigma_i\sigma_{i+1}$, for $i=1,\ldots, n-2$,
\item[(4)] $\rho_{i+1}\sigma_i\sigma_{i+1} = \sigma_i\sigma_{i+1}\rho_i$, for $i=1,\ldots, n-2$,
\item[(5)] $\sigma_i\rho_i = \rho_i \sigma_i$, for $i=1,\ldots, n-1$.
\end{itemize}

\begin{Obs}
We note that we defined here $GVT_n$, the ``twin version'' for Gauss virtual braids $GVB_n$, and related them as one being the quotient of the another. We also may see that the flat welded braid group $FWB_n$ is a quotient of the welded twin group $WT_n$.

It is also an interesting fact, to observe in the diagram given in Figure~\ref{quocientes}, that the virtual braid and virtual twin groups have the flat virtual braid group $FVB_n$ and the Gauss virtual braid group $GVB_n$ as a common quotient group.  
\end{Obs}

Now we move to the notion of crystallographic groups, that will be very important in this paper. 
First, we give the following characterization of it (see \cite[Lemma~8]{GGO1} and also \cite[Theorem~2.1.4]{DEKIMPE1996}). 

\begin{Def}\label{lema8}
	A group $\Pi$ is said to be a crystallographic group if there is an integer $n\in \N$ and a short exact sequence
	\begin{equation}\label{sec}
		0\longrightarrow\Z^n\longrightarrow\Pi\stackrel{\zeta}{\longrightarrow}\Phi\longrightarrow 1
	\end{equation}
	such that:
	\begin{enumerate}
		\item[(i)] $\Phi$ is finite.
		\item[(ii)] The integral representation $\Theta\colon\Phi\longrightarrow \aut(\Z^n)$, induced by conjugation on $\Z^n$ and defined by $\Theta(\varphi)(x)=\pi x \pi^{-1}$, where $x\in \Z^n$, $\varphi\in\Phi$ and $\pi\in\Pi$ is such that $\zeta(\pi)=\varphi$, is faithful. 
	\end{enumerate}
The group $\Phi$ is the holonomy group of $\Pi$ and $n$ is the dimension of $\Pi$.
\end{Def}

Crystallographic groups play an important role in the study of the groups of isometries of Euclidean spaces. 
Torsion free crystallographic groups are called \textit{Bieberbach groups}, they are of very importance in Riemannian geometry, because there is a correspondence between the class of Bieberbach groups and the class of compact flat Riemannian manifolds, through the fundamental group (see \cite{22} and \cite{Wolf2011}).  

As far as we know, the first connection between crystallographic groups and braid groups was established in \cite{GGO1} for the case of classical Artin braid groups, and several properties of it were studied. 
In addition to showing that the quotient group $B_n/\Gamma_2(P_n)$ of the Artin braid group $B_n$ by the commutator subgroup of its pure Artin braid group $P_n$ is a crystallographic group, in \cite{GGO1} it was studied the torsion, the conjugacy classes of its finite-order elements, and the realization of abstract finite groups as subgroups of $B_n/\Gamma_2(P_n)$. 
Later, in \cite{gamma2}, it was studied the conjugacy problem and the realization of virtually cyclic groups in the group $B_n/\Gamma_2(P_n)$.
The surprising connection between Artin braid groups and crystallographic was extended to generalized braid groups associated to arbitrary complex reflection groups in \cite{M2016} (see also \cite{BM2020}), to orientable closed surfaces of positive genus in \cite{GGOP2021} and to punctured surfaces (orientable or not) in \cite{D2020}. 
It was proved in \cite[Proposition~17]{GGOP2021} that if $M$ is the sphere or a closed non-orientable surface then,  for every $n\geq 1$, the quotient $B_n(M)/\Gamma_2(P_n(M))$ is not a crystallographic group. 

Hence, it is natural to ask which \textit{braid-like groups} have subquotients isomorphic to crystallographic groups. 
In this paper, we deal with this problem for the case of virtual braid and virtual twin groups and related groups. 
We shall prove that if $G$ is one of the groups $VB_n$, $VT_n$, $WB_n$, $FVB_n$ and $FWB_n$ and if $N$ is the respective pure subgroup in each case then $G/\Gamma_2(N)$ is a crystallographic group for every $n\geq 1$. 
In each case, we explore structural properties of the respective crystallographic quotient as described below. 
We  note that in \cite[Section~3]{Betal2021} the authors proved that there is a decomposition $VB_n/\Gamma_2(VP_n)=VP_n/\Gamma_2(VP_n)\rtimes S_n$ and that it is a crystallographic group, and also they studied linear representations of this group. 
In this paper we reprove the connection between virtual braid groups and crystallographic groups and then we deeply study structural aspects of $VB_n/\Gamma_2(VP_n)$.

This paper is organized as follows. 
Section~\ref{new} is dedicated to study the existence of finite order elements and the conjugacy classes of elements of a group $G$ under the following conditions.
We consider a short exact sequence $\xymatrix{
		1 \ar[r] & \Z X \ar[r]^-{} &  G \ar[r]^-{{\pi}} &  S_n \ar[r] & 1}$
where $\Z X$ is a free Abelian group on a finite set $X$ and $S_n$ acts on the free Abelian group $\Z X$ as a faithful permutation representation.
In Theorem~\ref{teo2_geral} we characterize finite order elements in $G$. The study of conjugacy classes of elements in $G$ are given in
Theorem~\ref{teo3_geral} and Theorem~\ref{teo4_geral}.
To prove these three results, we use equations in the free Abelian group $\Z X$ in order to transform our initial question into a problem of finding a solution of a system of equations over the integers.
We use these general theorems in the rest of the paper when $G$ is one of the quotients $VB_n/\Gamma_2(VP_n)$ or $VT_n/\Gamma_2(PVT_n)$.
We note that we may apply the results of Section~\ref{new} for other groups, for instance for the crystallographic groups considered in  \cite{D2020}, \cite{GGO1} and \cite{GGOP2021}.

In Section~\ref{groupvbn} we recall first the presentation of (pure) virtual braid groups, and then we prove in Theorem~\ref{teo1} that $VB_n/\Gamma_2(VP_n)$ is a crystallographic group. Next, we show in Proposition~\ref{prop:embedofvbn} that $VB_n/\Gamma_2(VP_n)$ contains a copy of the quotient of the classical Artin braid group $B_n/\Gamma_2(P_n)$.
The rest of Subsection~\ref{subsec:estrutura} is devoted to exhibit some structural properties of the quotient $VB_n/\Gamma_2(VP_n)$ like classification of finite order elements (Theorem~\ref{teo2}), the study of conjugacy classes of elements (Theorems~\ref{teo3} and~\ref{teo4}) and the realization of a family of infinite virtually cyclic groups (Theorem~\ref{teo5}). 
In Subsection~\ref{subsec:realizaBieberbach} we realize Bieberbach groups in $VB_n/\Gamma_2(VP_n)$. Then, in Theorem~\ref{anosov} we consider the compact flat Riemannian manifolds that have fundamental group a Bieberbach group of Theorem~\ref{teo6}, we analyze its holonomy representation and determine in which cases the aforementioned manifolds are orientable, admits Anosov diffeomorphisms and K\"ahler structure.
In Subsection~\ref{subsec:kbn} we consider another subgroup of $VB_n$ denoted by $KB_n$ and conclude that $VB_n/\Gamma_2(KB_n)$ is not a crystallographic group.

Virtual twin groups modulo the commutator subgroup of pure virtual twin groups are considered in Section~\ref{groupvtn}.  
In this section we state and prove many results similar to the ones obtained for the virtual braid group, in particular in Theorem~\ref{teo7} is proved that $VT_n/\Gamma_2(PVT_n)$ is a crystallographic group. 
The rest of the section is dedicated to analyze finite-order elements, conjugacy classes of elements and realization of some infinite virtually cyclic groups in $VT_n/\Gamma_2(PVT_n)$.

The conjugacy problem, formulated by Max Dehn, is one of the fundamental decision problems in group theory. 
We highlight that, from Theorems~\ref{teo3}~and~\ref{teo4} (resp. Theorems~\ref{teo9}~and~\ref{teo10}), we give necessary and sufficient conditions for a pair of words in $VB_n/\Gamma_2(P_n)$ (resp. $VT_n/\Gamma_2(PVT_n)$) to be conjugate. So, we may use these results
to find elements that are not conjugate in $VB_n$ (resp. $VT_n$). 
In this way, we obtain an approach to solve the conjugacy problem in the virtual braid group $VB_n$ (resp. virtual twin group $VT_n$).

Finally, in Section~\ref{quotgroupvbn}, we consider other braid-like groups modulo the respective commutator subgroup in each case. 
In Subsection~\ref{subsec:wf} 
we show in Theorem~\ref{teo:isowbn/gamma} (resp. in Theorem~\ref{teo:isofvbn/gamma}) that the corresponding quotients of the welded and unrestricted virtual braid group (resp. flat virtual and flat welded braid group) are isomorphic to $VB_n/\Gamma_2(VP_n)$ (resp. $VT_n/\Gamma_2(PVT_n)$). 
%
%
%
In Subsection~\ref{subsec:ext} we consider a group that is not a quotient of the virtual braid group, neither of the virtual twin group. It is called the \textit{pure extended loop braid group}, and we show in Theorem~\ref{teo1ext} that it has a subquotient that is crystallographic.




\subsection*{Acknowledgements}

The authors would like to thank Paolo Bellingeri and John Guaschi for helpful comments, in particular for pointing out the isomorphisms described in Theorem~\ref{teo:isowbn/gamma}, avoiding unnecessary computations in the quotients of  welded and unrestricted braid groups. This work is part of the PhD thesis of the second author that was supported by CAPES.
We would like to thank an anonymous referee for a thorough reading and who generously made a number of suggestions throughout the text, in particular for proposing us to write the general results of Section~\ref{new}.


\section{Finite order elements and conjugacy classes in free abelian-by-finite groups}\label{new}


Let $G$ be a group that acts on a set $X$. 
By an abuse of language, we just said that $g\in G$ acts on $X$ meaning that the group generated by $ g $ acts on $X$.
Let $\Z X$ be the free Abelian group on the finite set X. 
Consider the following short exact sequence
\begin{equation*}
	\xymatrix{
		1 \ar[r] & \Z X \ar[r]^-{} &  G \ar[r]^-{{\pi}} &  F \ar[r] & 1}
\end{equation*}
where $F$ is a finite group that acts on the free Abelian group $\Z X$ as a permutation representation (i.e. $F$ acts permuting the elements of the basis of $\Z X$) that is  faithful.
We will denote a transversal of the action of $g$ on $X$ by $T_{\theta}$ and the orbit of an element $\lambda_{i}\in X$ by $\mathcal{O}_{\theta}(\lambda_{i})$, where $g\lambda_ig^{-1}=\lambda_{\theta(i)}$ with $\pi(g)^{-1}=\theta$.

Our first result in this section characterize finite order elements in $G$.

\begin{Teo}\label{teo2_geral}
Let $\Z X$ be the free Abelian group on the finite set $X$. 
Consider the following short exact sequence
\begin{equation}\label{sqpura}
	\xymatrix{
		1 \ar[r] & \Z X \ar[r]^-{} &  G \ar[r]^-{{\pi}} &  S_n \ar[r] & 1}
\end{equation}
where $S_n$ acts on the free Abelian group $\Z X$ as a faithful permutation representation. 
Let $\alpha\in G$ of order $t$. 
Let $\pi(\alpha)^{-1}=\theta$ and let $T_{\theta}$ be a transversal of the action of $\alpha$ on $\{\lambda_i\mid \lambda_i\in X\}$. The element $\prod_{\lambda_i\in X}\lambda_i^{a_i}\alpha$, where $a_{i}\in\Z$, has order $t$ in $G$ if and only if $\sum_{\lambda_{r}\in \mathcal{O}_{\theta}(\lambda_{i_{j}})}a_{r}=0$ for all $\lambda_{i_{j}}\in T_{\theta}$.
\end{Teo}

\begin{proof}
We will analyze under which conditions of $a_{i}\in\Z$ the element $\prod_{\lambda_i\in X}\lambda_i^{a_i}\alpha$ has order $t$ in $G$.
Consider $\pi(\alpha)^{-1}=\theta$ and let $T_{\theta}$ be a transversal of the action of $\alpha$ on $\{\lambda_{i}\mid \lambda_i\in X\}$. Hence,
\begin{align*}
\pare\prod_{\lambda_i\in X}\lambda_i^{a_i}\alpha\pard^t&=\prod_{\lambda_i\in X}\lambda_i^{a_i}\alpha\prod_{\lambda_i\in X}\lambda_i^{a_i}\alpha\prod_{\lambda_i\in X}\lambda_i^{a_i}\alpha\cdots\prod_{\lambda_i\in X}\lambda_i^{a_i}\alpha\\
&=\prod_{\lambda_i\in X}\lambda_i^{a_i}\alpha\prod_{\lambda_i\in X}\lambda_i^{a_i}\alpha^{-1}\alpha^2\prod_{\lambda_i\in X}\lambda_i^{a_i}\alpha^{-2}\alpha^{3}\cdots\alpha^{t-1}\prod_{\lambda_i\in X}\lambda_i^{a_i}\alpha^{-(t-1)}\alpha^t\\
&=\prod_{\lambda_i\in X}\lambda_i^{a_i}\prod_{\lambda_i\in X}\lambda_{\theta(i)}^{a_{i}}\prod_{\lambda_i\in X}\lambda_{\theta^2(i)}^{a_{i}}
\cdots\prod_{\lambda_i\in X}\lambda_{\theta^{t-1}(i)}^{a_{i}}.
\end{align*}

Since $\{\lambda_{i}\mid \lambda_i\in X\}$ is a basis of the free Abelian group $\Z X$ we note that, choosing the transversal $T_{\theta}=\{\lambda_{i_1},\lambda_{i_2},\dots,\lambda_{i_m}\}$, to solve the following equation in $\Z X$
\begin{align*}
\pare\prod_{\lambda_i\in X}\lambda_i^{a_i}\alpha\pard^t&=\prod_{\lambda_i\in X}\lambda_i^{a_i}\prod_{\lambda_i\in X}\lambda_{\theta(i)}^{a_{i}}\prod_{\lambda_i\in X}\lambda_{\theta^2(i)}^{a_{i}}
\cdots\prod_{\lambda_i\in X}\lambda_{\theta^{t-1}(i)}^{a_{i}}=1
\end{align*}
is equivalent to solve the system of integer equations
$$
\begin{cases}
	a_{i_1}+a_{\theta^{-1}(i_1)}+\dots+a_{\theta^{-(t-1)}(i_1)}=0&\\
	a_{i_2}+a_{\theta^{-1}(i_2)}+\dots+a_{\theta^{-(t-1)}(i_2)}=0&\\
	\vdots\quad\quad\vdots\quad\quad\vdots\quad\quad\vdots\quad\quad\vdots\quad\quad\vdots\quad\quad&\\
	a_{i_m}+a_{\theta^{-1}(i_r)}+\dots+a_{\theta^{-(t-1)}(i_r)}=0&
	\end{cases}
$$
i.e. $\pare\prod_{\lambda_i\in X}\lambda_i^{a_i}\alpha\pard^t=1$ if and only if $\sum_{\lambda_{r}\in \mathcal{O}_{\theta}(\lambda_{i_j})}a_{r}=0$ for all $\lambda_{i_j}\in T_{\theta}$.
\end{proof}

We register as a corollary the case in which the group $G=\Z X \rtimes S_n$ of Theorem~\ref{teo2_geral} is a semi-direct product. In this way we characterize all finite order elements in $G$.

\begin{Cor}\label{cor_teo2geral}
Let $\Z X$ be the free Abelian group on the finite set $X$. 
Consider the following split short exact sequence
\begin{equation*}
	\xymatrix{
		1 \ar[r] & \Z X \ar[r]^-{} &  G \ar[r]^-{{\pi}} &  S_n \ar[r] & 1}
\end{equation*}
where $S_n$ acts on the free Abelian group $\Z X$ as a faithful permutation representation. 
Let $\alpha\in S_n$ of order $t$. 
Let $\pi(\alpha)^{-1}=\theta$ and let $T_{\theta}$ be a transversal of the action of $\alpha$ on $\{\lambda_i\mid \lambda_i\in X\}$. The element $\prod_{\lambda_i\in X}\lambda_i^{a_i}\alpha$, where $a_{i}\in\Z$, has order $t$ in $G$ if and only if $\sum_{\lambda_{r}\in \mathcal{O}_{\theta}(\lambda_{i_{j}})}a_{r}=0$ for all $\lambda_{i_{j}}\in T_{\theta}$.
\end{Cor}


Let $\alpha\in G$. A necessary condition for a given coset of an element $y\in G$ to belong to the conjugacy class of $\alpha$ is that the permutations $\pi(y)$ and $\pi(\alpha)$ need to have the same cycle type, where $\pi$ is the surjective homomorphism given in (\ref{sqpura}).
Now, consider $y\in G$ such that $\pi(y)$ and $\pi(\alpha)$ have the same cycle type. So, there is $\hat{\gamma}\in S_n$ such that $\hat{\gamma}\pi(y)\hat{\gamma}^{-1}=\pi(\alpha)$. 
Since $\pi$ is surjective, then there is $\gamma \in G$ such that $\tilde{y}=\gamma y \gamma^{-1}$ belongs to $\pi^{-1}(\pi(\alpha))$.
We note that $\tilde{y}$ is not necessarily conjugate to $\alpha$. For an example of this situation see \cite[Remark~2.1]{gamma2}. 
Our first result about the conjugacy problem in $G$ is as follows.

\begin{Teo}\label{teo3_geral}
Let $\Z X$ be the free Abelian group on the finite set $X$. 
Consider the following short exact sequence
\begin{equation*}
	\xymatrix{
		1 \ar[r] & \Z X \ar[r]^-{} &  G \ar[r]^-{{\pi}} &  S_n \ar[r] & 1}
\end{equation*}
where $S_n$ is a finite group that acts on the free Abelian group $\Z X$ as a faithful permutation representation.
Let $\alpha\in G$, let $\theta=\pi(\alpha)^{-1}$ and let $T_{\theta}$ be a transversal of the action of $\alpha$ on $\{\lambda_{i}\mid \lambda_i\in X\}$. 
The elements $\prod_{\lambda_i\in X}\lambda_{i}^{a_{i}}\alpha$ and $\prod_{\lambda_{i_j}\in T_{\theta}}\lambda_{i_j}^{S_{i_j}}\alpha$ are conjugate in $G$, where $S_{i_j}=\sum_{\lambda_{r}\in\mathcal{O}_{\theta}(\lambda_{i_{j}})}{a_{r}}$ for all $\lambda_{i_j}\in T_{\theta}$.
\end{Teo}

\begin{proof}
In order to prove that $\prod_{\lambda_i\in X}\lambda_{i}^{a_{i}}\alpha$ and $\prod_{\lambda_{i_j}\in T_{\theta}}\lambda_{i_j}^{S_{i_j}}\alpha$ are conjugate in $G$ we will show that there is $Y\in \Z X$ such that
$$
Y\prod_{\lambda_i\in X}\lambda_{i}^{a_{i}}\alpha Y^{-1}=\prod_{\lambda_{i_j}\in T_{\theta}}\lambda_{i_j}^{S_{i_j}}\alpha,
$$ 
or equivalently,  
$$
Y\prod_{\lambda_i\in X}\lambda_{i}^{a_{i}}\alpha Y^{-1}\alpha^{-1}=\prod_{\lambda_{i_j}\in T_{\theta}}\lambda_{i_j}^{S_{i_j}}.
$$ 

Using the basis of the free Abelian group $\Z X$ we consider $Y=\prod_{\lambda_i\in X}\lambda_{i}^{y_{i}}\in \Z X$, where $y_{i}\in\Z$. 
We will analyze conditions on the elements $y_{i}\in\Z$ in such a way $\prod_{\lambda_i\in X}\lambda_{i}^{a_{i}}\alpha$ and $\prod_{\lambda_{i_j}\in T_{\theta}}\lambda_{i_j}^{S_{i_j}}\alpha$ are conjugate in $G$. First, we note that
\begin{align*}
&Y\prod_{\lambda_i\in X}\lambda_{i}^{a_{i}}\alpha Y^{-1}\alpha^{-1}=\prod_{\lambda_{i_j}\in T_{\theta}}\lambda_{i_j}^{S_{i_j}}\\
&\Leftrightarrow Y\prod_{\lambda_i\in X}\lambda_{i}^{a_{i}}\prod_{\lambda_i\in X}\lambda_{\theta(i)}^{-y_{i}}=\prod_{\lambda_{i_j}\in T_{\theta}}\lambda_{i_j}^{S_{i_j}},
\end{align*}
and since the elements $\lambda_{i}$ are generators of the free Abelian group $\Z X$ then we conclude that
$$
Y\prod_{\lambda_i\in X}\lambda_{i}^{a_{i}}\alpha Y^{-1}\alpha^{-1}=\prod_{\lambda_{i_j}\in T_{\theta}}\lambda_{i_j}^{S_{i_j}}
$$
\begin{align}\label{sistconj}
\Leftrightarrow\begin{cases}
	y_{i}+a_{i}-y_{\theta^{-1}(i)}=0,&\mbox{ if }\lambda_{i}\notin T_{\theta},\\
	y_{i_j}+a_{i_j}-y_{\theta^{-1}(i_j)}=S_{i_j},&\mbox{ if }\lambda_{i_j}\in T_{\theta}.
\end{cases}
\end{align}

The system of equations~(\ref{sistconj}) has a solution if and only if the following subsystems of equations has a solution   
	\begin{align}\label{subsistconj}
	\Leftrightarrow\begin{cases}
		y_{\theta^t(i_j)}+a_{\theta^t(i_j)}-y_{\theta^{-(1+m-t)}(i_j)}=0,\\
		y_{i_j}+a_{i_j}-y_{\theta^{-1}(i_j)}=S_{i_j}.
	\end{cases}
\end{align}
for all $1\le t<m$ and for all $\lambda_{i_j}\in T_{\theta}$ where $m=|\mathcal{O}_{\theta}(\lambda_{i_j})|$.

So, in order to prove the statement of this theorem, we will show that the subsystems of equations~(\ref{subsistconj}) admits solution. Let $\lambda_{i_j}\in T_{\theta}$ and consider $e_1=\lambda_{i_j},e_2=\lambda_{\theta(i_j)},\dots,e_m=\lambda_{\theta^{m-1}(i_j)}$ where $m=|\mathcal{O}_{\theta}(\lambda_{i_j})|$. Thus, the permutation matrix of the orbit of $\lambda_{i_j}$ by the action of $\theta$ with respect to the ordered set $\{e_1,e_2,\dots,e_m\}$ is given by the matrix $M_{i_j}$
$$M_{i_j}=
\begin{pmatrix}
	0&0&\ldots&0&1\\
	1&0&\ldots&0&0\\
	0&1&\ldots&0&0\\
	\vdots&\vdots&\ddots&\vdots&\vdots\\ 
	0&0&\ldots&1&0
\end{pmatrix}.
$$

So, for every $\lambda_{i_j}\in T_{\theta}$, we can write the subsystems of equations~(\ref{subsistconj}) as follows, 
$$
(\star)\quad [I_{i_j}-M_{i_j}]
\begin{pmatrix}
	y_{i_j}\\
	y_{\theta(i_j)}\\
	\vdots\\
	y_{\theta^{m-1}(i_j)}
\end{pmatrix}=
\begin{pmatrix}
	-a_{i_j}+S_{i_j}\\
	-a_{\theta(i_j)}\\
	\vdots\\
	-a_{\theta^{m-1}(i_j)}
\end{pmatrix}
$$ 
where $I_{i_j}=(Id)_{m\times m}$ and $m=|\mathcal{O}_{\theta}(\lambda_{i_j})|$. Analyzing the coefficient matrix and the extended matrix associated with this system, we conclude that $(\star)$ have a solution if and only if $S_{i_j}=\sum_{\lambda_{r}\in\mathcal{O}_{\theta}(\lambda_{i_j})}{a_{r}}$ for all $\lambda_{i_j}\in T_{\theta}$. 
\end{proof}


With the notations of Theorem~\ref{teo3_geral} let $T'_{\theta}$ be another transversal of the action of $\langle\theta\rangle$ on the set $X$. Suppose that $\lambda_{i_j'}\in T'_{\theta}$ is the representative of the orbit of $\lambda_{i_j}\in T_{\theta}$. 
From Theorem~\ref{teo3_geral}, if $S_{i_j}=\sum_{\lambda_{i_j}\in\mathcal{O}_{\theta}(\lambda_{i_j})}{a_{i_j}}=\sum_{\lambda_{i_j}\in\mathcal{O}_{\theta}(\lambda_{i_j'})}{a_{i_j}}$ then $\prod_{\lambda_{i_j}\in T_{\theta}}{\lambda_{i_j}^{S_{i_j}}} \alpha$ and $\prod_{\lambda_{i_j'}\in T'_{\theta}}{\lambda_{i_j'}^{S_{i_j}}} \alpha$ are conjugate.
Let $\alpha\in G$ satisfying the conditions of Theorem~\ref{teo3_geral}. 
Let $\pi(\alpha^{-1})=\theta$, we choose a transversal $T_{\theta}$ of the action of $\langle \theta \rangle$ on the set $X$. 
In the next theorem we study which kind of elements belongs to the conjugacy class of $ \prod_{\lambda_{i_j}\in T_{\theta}}  {\lambda_{i_j}^{ S_{i_j} } }\alpha$, having the same description using the chosen transversal $T_{\theta}$, where $S_{i_j}=\sum_{\lambda_{r}\in\mathcal{O}_{\theta}(\lambda_{i_{j}})}{a_{r}}$ for all $\lambda_{i_j}\in T_{\theta}$. 
We will denote by $C_{F}(\theta)$ the centralizer of $\theta$ in $F$.

\begin{Teo}\label{teo4_geral}
Let $\Z X$ be the free Abelian group on the finite set $X$. 
Consider the following short exact sequence
\begin{equation*}
	\xymatrix{
		1 \ar[r] & \Z X \ar[r]^-{} &  G \ar[r]^-{{\pi}} &  S_n \ar[r] & 1}
\end{equation*}
where $S_n$ acts on the free Abelian group $\Z X$ as a faithful permutation representation.
Let $\alpha\in G$, let $\theta=\pi(\alpha)^{-1}$ and let $c\in C_{F}(\theta)$. 
Consider a transversal $T_{\theta}$ of the action of $\alpha$ on the set of generators $\{\lambda_{i}\mid \lambda_i\in X\}$, $\tilde{c}\in G$ such that $\pi(\tilde{c})^{-1}=c$ and $\tilde{c}\alpha\tilde{c}^{-1}=\prod_{\lambda_i\in X}\lambda_i^{c_i}\alpha$. 
Then, the elements $\prod_{\lambda_{r}\in T_{\theta}}\lambda_{r}^{a_{r}}\alpha$ and $\prod_{\lambda_{r}\in T_{\theta}}\lambda_{r}^{b_{r}}\alpha$ are conjugate in $G$ if and only if

\begin{enumerate}
\item[$(i)$] $\sum_{\lambda_i\in\mathcal{O}(\lambda_r)}c_i+a_{r}=b_{r}$, for $\lambda_{c^{-1}(r)}\in\mathcal{O}_{\theta}(\lambda_{r})$ where $\lambda_{r}\in T_{\theta}$.

\item[$(ii)$] $\sum_{\lambda_i\in\mathcal{O}(\lambda_r)}c_i+a_{c^{-1}(r)}=b_{r}$, for $\lambda_{c^{-1}(r)}\not\in\mathcal{O}_{\theta}(\lambda_{r})$ where  $a_{c^{-1}}=\sum_{\lambda_{i}\in\mathcal{O}_{\theta}(\lambda_{c^{-1}(r)})}a_{i}$ and $\lambda_{r}\in T_{\theta}$.
\end{enumerate}	
\end{Teo}

\begin{proof}
To verify that $\prod_{\lambda_{r}\in T_{\theta}}\lambda_{r}^{a_{r}}\alpha$ and $\prod_{\lambda_{r}\in T_{\theta}}\lambda_{r}^{b_{r}}\alpha$ are conjugate in $G$ we need to find an element $Y\in \Z X$ such that 
\begin{equation}\label{siscentra}
Y\tilde{c}\prod_{\lambda_{r}\in T_{\theta}}\lambda_{r}^{a_{r}}\alpha\tilde{c}^{-1}Y^{-1}=\prod_{\lambda_{r}\in T_{\theta}}\lambda_{r}^{b_{r}}\alpha,
\end{equation}
where $\tilde{c}\in G$ is such that $\pi(\tilde{c})^{-1}=c$.

Consider $Y=\prod_{\lambda_i\in X}{\lambda_{i}^{y_{i}}}\in \Z X$. We will analyze conditions on $Y$ so that equation~(\ref{siscentra}) holds. To do that we use the following equivalences of equations in the free Abelian group $\Z X$ in order to transform our initial question into a problem of finding a solution of a system of equations over the integers,
\begin{align*}
&Y\tilde{c}\prod_{\lambda_{r}\in T_{\theta}}\lambda_{r}^{a_{r}}\alpha\tilde{c}^{-1}Y^{-1}=\prod_{\lambda_{r}\in T_{\theta}}\lambda_{r}^{b_{r}}\alpha\\
\Leftrightarrow & Y\tilde{c}\prod_{\lambda_{r}\in T_{\theta}}\lambda_{r}^{a_{r}}\alpha\tilde{c}^{-1}Y^{-1}\alpha^{-1}=\prod_{\lambda_{r}\in T_{\theta}}\lambda_{r}^{b_{r}}\\
\Leftrightarrow & Y\prod_{\lambda_{r}\in T_{\theta}}\lambda_{c(r)}^{a_{r}}\tilde{c}\alpha\tilde{c}^{-1}Y^{-1}\alpha^{-1}=\prod_{\lambda_{r}\in T_{\theta}}\lambda_{r}^{b_{r}},
\end{align*}
and since $\pi(\tilde{c})^{-1}=c$ and $c\in C_{F}(\theta)$ we have $\tilde{c}\alpha\tilde{c}^{-1}=\prod_{\lambda_i\in X}\lambda_i^{c_i}\alpha$
\begin{align*}
\Leftrightarrow &Y\prod_{\lambda_{r}\in T_{\theta}}\lambda_{c(r)}^{a_{r}}\prod_{\lambda_i\in X}\lambda_i^{c_i}\alpha Y^{-1}\alpha^{-1}=\prod_{\lambda_{r}\in T_{\theta}}\lambda_{r}^{b_{r}}\\	
\Leftrightarrow &\prod_{\lambda_i\in X}{\lambda_{i}^{y_{i}}}\prod_{\lambda_i\in X}\lambda_i^{c_i}\prod_{\lambda_{r}\in T_{\theta}}\lambda_{c(r)}^{a_{r}}\prod_{\lambda_i\in X}{\lambda_{\theta(i)}^{-y_{i}}}=\prod_{\lambda_{r}\in T_{\theta}}\lambda_{r}^{b_{r}}.
\end{align*}

Suppose that $\lambda_{c^{-1}(r)}\notin\mathcal{O}_{\theta}(\lambda_{r})$. 
We assert that $a_{c^{-1}(r)}$, the exponent of $\lambda_{c^{-1}(r)}$, is equal to the exponent of some element in the transversal $T_{\theta}$. In fact, if $\lambda_{c^{-1}(r)}\notin T_{\theta}$, then there is $\lambda_{r_1}\in T_{\theta}$ such that $\lambda_{c^{-1}(r)}\in\mathcal{O}_{\theta}(\lambda_{r_1})$. From item $(ii)$ the equality $a_{c^{-1}(r)}=\sum_{\lambda_{i}\in\mathcal{O}_{\theta}(\lambda_{c^{-1}(r)})}a_{i,j}$ holds and since $\lambda_{c^{-1}(r)}\in\mathcal{O}_{\theta}(\lambda_{r_1})$, we have $a_{c^{-1}(r)}=a_{r_1}$. So, the last equation is true if the following system of equations has a solution 
\begin{align}\label{subsistcenconj}
\begin{cases}
y_{\theta^k(r)}+c_{\theta^k(r)}-y_{\theta^{-(1+m-k)}(r)}=0,\\
y_{r}+c_{r}+a_{c^{-1}(r)}+-y_{\theta^{-1}(r)}-b_{r}=0.	
\end{cases} 
\end{align}
for $1\le k<m$, where $\lambda_{r}\in T_{\theta}$ and $m=|\mathcal{O}_{\theta}(\lambda_{r})|$. Applying the same argument given in the proof of Theorem~\ref{teo3_geral}, we may conclude that equation~(\ref{subsistcenconj}) is equivalent to 
$$
(I)\quad [I_{r,s}-M_{r,s}]
\begin{pmatrix}
	x_{r}\\
	x_{\theta(r)}\\
	\vdots\\
	x_{\theta^{m-1}(r)}
\end{pmatrix}=
\begin{pmatrix}
	-c_r-a_{c^{-1}(r)}+b_{r}\\
	-c_\theta(r)\\
	\vdots\\
	-c_\theta^{m-1}(r)
\end{pmatrix}
$$
where $m=|\mathcal{O}_{\theta}(\lambda_{r})|$, $I_{r}=(Id)_{m\times m}$ and $M_{r}$ is the permutation matrix of the orbit of $\lambda_{r}$ by the action of $\theta$ on $\{\lambda_{i}\mid \lambda_i\in X\}$ with respect to the ordered set $\{e_1=\lambda_{r},e_2=\lambda_{\theta(r)},\dots,e_m=\lambda_{\theta^{m-1}(r)}\}$. Analyzing the coefficient matrix and the extended matrix associated to this system, we concluded that $(I)$ have solution if and only if $\sum_{\lambda_i\in\mathcal{O}(\lambda_r)}c_i+a_{c^{-1}(r)}=b_{r}$ proving item $(ii)$. 
Similarly, if $\lambda_{c^{-1}(r)}\in\mathcal{O}_{\theta}(\lambda_{r})$ for $\lambda_{r}\in T_{\theta}$ the system $(I)$ have solution if and only if $\sum_{\lambda_i\in\mathcal{O}(\lambda_r)}c_i+a_{r}=b_{r}$.  
\end{proof}

We note that the proofs of Theorems~\ref{teo3_geral} and~\ref{teo4_geral} are analogous, respectively, to the proofs of Theorems~$1.1$ and~$1.2$ formulated in \cite{gamma2} for a similar quotient of the Artin braid group $B_n$.

To finish this section, we describe the steps to compare whether two given elements are conjugate in $G$. 
Let $\beta\in G$. Let $\alpha$ be an arbitrary element different from $\beta$, we would like to decide whether $\alpha$ is conjugate to $\beta$ in $G$. 
Follows from the existence of the short exact sequence given in equation~(\ref{sqpura}) that a necessary condition to $\alpha$ and $\beta$ to be conjugate is that their permutations related have the same cycle type.
\begin{itemize}

	\item[Step 1.] Since $\alpha$ and $\beta$ have permutations with the same cycle type and $\pi$ is surjective then $\alpha$ is conjugate to an element $\alpha_1$ in $\pi^{-1}(\pi(\beta))$.
 
\item[Step 2.]  
Suppose that $\alpha_1=\prod_{\lambda_i \in X}{\lambda_{i}^{a_{i}}} \beta\in \pi^{-1}(\pi(\beta))$, where $a_{i}\in\Z$ for all $\lambda_i \in X$. 
Let $\theta=\pi(\beta^{-1})$. From Theorem~\ref{teo3_geral} $\alpha_1$ is conjugate to an element $\alpha_2=\prod_{\lambda_{i_j}\in T_{\theta}}{\lambda_{i_j}^{S_{i_j}}} \beta$ in $G$, where $T_{\theta}$ is a transversal of the action of $\langle \theta \rangle$ on the set $X$, $\mathcal{O}_{\theta}(\lambda_{i_j})$ is the orbit of the element $\lambda_{i_j}$ under this action and $S_{i_j}=\sum_{\lambda_{i}\in\mathcal{O}_{\theta}(\lambda_{i_j})}{a_{i}}$, for all $\lambda_{i_j}\in T_{\theta}$.

\item[Step 3.]
Finally, using Theorem~\ref{teo4_geral}, we are able to decide if $\alpha_2$ is conjugate to $\beta$, solving a system of equations over the integers.  
As part of this procedure we need to choose $\tilde{C}\in\pi^{-1}(C_{S_n}(\theta))$, where $C_{S_n}(\theta)$ is the centralizer of $\theta$ in $S_n$, and to write $\tilde{C}\beta\tilde{C}^{-1}$ in the form $\prod_{\lambda_i\in X}\lambda_i^{c_i}\beta$. 
\end{itemize}

An example of the utilization of this method may be found in \cite[Example~2.7]{gamma2}.

\section{Crystallographic groups and virtual braids}\label{groupvbn}

In this section, we establish the connection between virtual braid groups and crystallographic groups.  
First, we rewrite a presentation of the virtual braid group $VB_n$ that will be very useful in this work. It was formulated  in \cite[p.798]{V2001} and restated in \cite{BB2009}.

\begin{Teo}[{\cite[Theorem~4]{BB2009}}]\label{apvbn}
The group $VB_n$ admits the following group presentation:
\begin{itemize}
    \item Generators: $\sigma_i$, $\rho_i$ where $i=1,2,\dots,n-1.$
    \item Relations:
       \subitem $\sigma_i\sigma_{i+1}\sigma_{i}=\sigma_{i+1}\sigma_{i}\sigma_{i+1}$, $i=1,2,\dots,n-2;$
       \subitem $\sigma_{i}\sigma_j=\sigma_{j}\sigma_i$, $\mid i-j\mid\ge 2;$
       \subitem $\rho_i\rho_{i+1}\rho_{i}=\rho_{i+1}\rho_{i}\rho_{i+1}$, $i=1,2,\dots,n-2;$
       \subitem  $\rho_{i}\rho_j=\rho_{j}\rho_i$, $\mid i-j\mid\ge 2;$
 \subitem $\rho_i^2=1$, $i=1,2,\dots,n-1;$
 \subitem  $\sigma_{i}\rho_j=\rho_{j}\sigma_i$, $ \mid i-j\mid\ge 2;$
\subitem $\rho_i\rho_{i+1}\sigma_{i}=\sigma_{i+1}\rho_{i}\rho_{i+1}$, $i=1,2,\dots,n-2.$
\end{itemize}
\end{Teo}

Let $n\geq 2$ and, let $\pi_p\colon VB_n\longrightarrow S_n$ be the homomorphism defined by ${\pi}_p(\sigma_i)={\pi}_p(\rho_i)=\tau_i$ for $i=1,\dots,n-1$ and $\tau_i=(i,i+1)$. For convention, in the case $n=1$ this homomorphism is the trivial one. 
The \textit{pure virtual braid group} $VP_n$ is defined to be the kernel of $\pi_p$, from which we obtain the following short exact sequence
\begin{equation}\label{seqvbn0} 
	\xymatrix{
		1 \ar[r] & VP_n \ar[r]^-{} &  VB_n \ar[r]^-{{\pi}_p} &  S_n \ar[r] & 1}.
\end{equation}
As mentioned in the first paragraph of Section~3 of \cite{B2004} the virtual braid group admits a decomposition as semi-direct product $VB_n=VP_n\rtimes S_n$, the map $\iota\colon S_n\longrightarrow VB_n$ given by $\iota(\tau_i)=\rho_i$, for $i=1,\dots,n-1$, is naturally a section for $\pi_p$. 




We state here, without proof, the following presentation of the pure virtual braid group. 

\begin{Teo}[{\cite[Theorem~1]{B2004}}]\label{apvpn}
  The group $VP_n$ admits the following group presentation: 
  \begin{itemize}
\item Generators:
      \subitem $\lambda_{i,j}$ for $1\le i\neq j\le n.$
\item Relations:
     \subitem $\lambda_{i,j}\lambda_{k,l}=\lambda_{k,l}\lambda_{i,j}$ for $i,j,k,l$ distinct.
     \subitem $\lambda_{k,i}(\lambda_{k,j}\lambda_{i,j})=(\lambda_{i,j}\lambda_{k,j})\lambda_{k,i}$ for $i,j,k$ distinct.
	\end{itemize}
\end{Teo}

From Lemma~1 of~\cite{B2004}, using the presentation of the virtual pure braid group given in Theorem~\ref{apvpn}, we conclude that the action  by conjugation of $S_n=\langle \rho_1,\ldots, \rho_{n-1} \rangle$ on $VP_n$ is described as follows: let $w\in S_n$ and let $\lambda_{i,j}$ be a generator of $VP_n$, with $1\le i\neq j\le n$, then 
\begin{equation}\label{conjvbn}
	w\cdot\lambda_{i,j}=w\lambda_{i,j}w^{-1}=\lambda_{w^{-1}(i),w^{-1}(j)}.
\end{equation}

\subsection{The virtual braid group modulo the commutator subgroup of the pure virtual braid group}\label{subsec:estrutura}




Recall that the pure virtual braid group $VP_n$ is defined to be the kernel of $\pi_p\colon VB_n \to S_n$. The split short exact sequence (\ref{seqvbn0}) provides, for every $n\geq 2$, the following short exact sequence
\begin{equation}\label{seqvbn} 
	\xymatrix{
		1 \ar[r] & VP_n/\Gamma_2(VP_n) \ar[r]^-{} &  VB_n/\Gamma_2(VP_n) \ar[r]^-{\overline{\pi}_p} &   S_n \ar[r] & 1.}
\end{equation}
We notice that $\iota$, the section for $\pi_p$ described in the paragraph below equation~(\ref{seqvbn0}), induces a homomorphism $\overline{\iota}\colon S_n\longrightarrow VB_n/\Gamma_2(VP_n)$ which is a section for $\overline{\pi}_p$, so we obtain the decomposition $VB_n/\Gamma_2(VP_n)=VP_n/\Gamma_2(VP_n)\rtimes S_n$.

We will apply this to the action by conjugation in $VB_n/\Gamma_2(VP_n)$.
Let $\nu\in VB_n/\Gamma_2(VP_n)$. 
Since $VB_n/\Gamma_2(VP_n)=VP_n/\Gamma_2(VP_n)\rtimes S_n$, we have $\nu=Aw$ for some $A\in VP_n/\Gamma_2(VP_n)$ and for some $w\in S_n$. 
 The action of $\nu$ on $\{\lambda_{i,j}\mid 1\le i\neq j\le n\}$ corresponds to the action of the group generated by $\overline{\pi}_p(\nu)$ on $\{\lambda_{i,j}\mid 1\le i\neq j\le n\}$ (see~(\ref{conjvbn})).  
We will denote a transversal of the action of $\nu$ on $\{\lambda_{i,j}\mid 1\le i\neq j\le n\}$ by $T_{\theta}$ and the orbit of an element $\lambda_{i,j}$ by $\mathcal{O}_{\theta}(\lambda_{i,j})$ where $\theta=\overline{\pi}_p(\nu)$.


The first result in this subsection is that $VB_n/\Gamma_2(VP_n)$ is a crystallographic group of dimension $n(n-1)$ and whose holonomy group is $S_n$. 
We  recall that this result is known, see \cite[Section~3]{Betal2021}. 
For a sake of completeness, we give a proof here.

\begin{Teo}\label{teo1}
	Let $n\ge 2$.  There is a split short exact sequence:
	\begin{equation}\label{seqvbn2} 
	\xymatrix{
		1 \ar[r] & \Z^{n(n-1)} \ar[r]^-{} &  VB_n/\Gamma_2(VP_n) \ar[r]^-{\overline{\pi}_p} &   S_n \ar[r] & 1}
\end{equation}
	and the middle group is a crystallographic group.

\end{Teo}
\begin{proof}
From the presentation of $VP_n$, see Theorem~\ref{apvpn}, it is straightforward to show that $VP_n/\Gamma_2(VP_n)$ is isomorphic to the free Abelian group $\Z^{n(n-1)}$ and  generated by the set 
\begin{equation}\label{basis}
\{\lambda_{i,j}\mid 1\le i\neq j\le n\}
\end{equation}

Thus, we get from the short exact sequence~(\ref{seqvbn})  the split short exact sequence~(\ref{seqvbn2}). 
Hence, according to Definition~\ref{lema8} applied to the short exact sequence~(\ref{seqvbn2}), rest to prove that the action of $S_n$ on $VP_n/\Gamma_2(VP_n)$ is injective to conclude that $VB_n/\Gamma_2(VP_n)$ is a crystallographic group.
The action of $S_n$ on $VP_n/\Gamma_2(VP_n)$ is induced from the action of $S_n$ on $VP_n$ given in (\ref{conjvbn}).

Let $\beta\in S_n$ such that $\beta\cdot\lambda_{i,j}=\lambda_{i,j}$ for all $\lambda_{i,j}\in VP_n/\Gamma_2(VP_n)$ , so we have  
\begin{align*}
	\beta\cdot\lambda_{i,j}=\lambda_{i,j}&\Leftrightarrow\beta\lambda_{i,j}\beta^{-1}=\lambda_{i,j}\\
	&\Leftrightarrow\lambda_{\beta^{-1}(i),\beta^{-1}(j)}=\lambda_{i,j},
\end{align*}
for all $1\le i\neq j\le n$. 
Since $\che \lambda_{i,j} \mid 1\le i\neq j\le n \chd$ is a basis of the free Abelian group $VP_n/\Gamma_2(VP_n)$, we get $\beta=1$. Therefore,
the action of $S_n$ on $VP_n/\Gamma_2(VP_n)$ is injective.
\end{proof}

The group $B_n$ naturally embeds in the virtual braid group $VB_n$, as mentioned in the proof of  Proposition~2 of~\cite{SW2017}. 
Below we show that  $B_n/\Gamma_2(P_n)$ is also embedded  in $VB_n/\Gamma_2(VP_n)$.


\begin{Prop}\label{prop:embedofvbn}
The embedding of $B_n$ into $VB_n$ induces an embedding of $B_n/\Gamma_2(P_n)$ into $VB_n/\Gamma_2(VP_n)$.  
\end{Prop}

\begin{proof}
Let $\varphi_n\colon B_n\hookrightarrow VB_n$ be the homomorphism defined in the proof of~Proposition~2 of~\cite{SW2017}. 
In~\cite[Proposition~4.8]{BEER2006} the authors showed that the homomorphism $\varphi_n$ restricts to a monomorphism $P_n\hookrightarrow VP_n$, given by
$$
A_{i,j}\mapsto(\lambda_{j-1,j}\ldots\lambda_{i+1,j})\lambda_{i,j}\lambda_{j,i}(\lambda_{j-1,j}\ldots\lambda_{i+1,j})^{-1}.
$$

The homomorphism $\varphi_n$ induces a homomorphism $\widetilde{{\varphi_n}}\colon B_n/\Gamma_2(P_n)\longrightarrow VB_n/\Gamma_2(VP_n)$. The restriction ${\varphi_n}_{\arrowvert_{P_n}}$ induces a monomorphism $\widetilde{{\varphi_n}}_{\arrowvert_{H}}\colon P_n/\Gamma_2(P_n)\longrightarrow VP_n/\Gamma_2(VP_n)$ where $H=P_n/\Gamma_2(P_n)$.
So, we obtain the following commutative diagram
\begin{eqnarray*}
\xymatrix{
1\ar[r] & P_{n}/\Gamma_2(P_n)\ar[d]^-{\widetilde{{\varphi_n}}_{\arrowvert_{H}}}\ar[r] & B_{n}/\Gamma_2(P_n)\ar[d]^-{\widetilde{{\varphi_n}}}\ar[r] & S_{n}\ar@{=}[d]\ar[r] & 1\\
1\ar[r]&VP_{n}/\Gamma_2(VP_n)\ar[r]&VB_{n}/\Gamma_2(VP_n)\ar[r]&S_{n}\ar[r]&1.}
\end{eqnarray*}
By using the five lemma we conclude that the homomorphism $\widetilde{{\varphi_n}}\colon B_n/\Gamma_2(P_n)\longrightarrow VB_n/\Gamma_2(VP_n)$ is injective.
\end{proof}

\begin{Obs}
It is natural to ask which elements in the quotient of the virtual braid group $VB_n/ \Gamma_2(VP_n)$ have finite order, as mentioned in \cite[Question~3.1]{Betal2021}. Also, in \cite[Question~3.2]{Betal2021} the authors asked about the realization of Bieberbach groups (i.e. torsion-free crystallographic groups) inside $VB_n/\Gamma_2(VP_n)$. 
It follows from Proposition~\ref{prop:embedofvbn} that in $VB_n/\Gamma_2(VP_n)$ there are infinitely many elements of finite order coming from $B_n/\Gamma_2(P_n)$ (see \cite[Theorem~3]{GGO1}) and also the Bieberbach groups with finite Abelian holonomy realized in $B_n/\Gamma_2(P_n)$ (see \cite{GGO1} and~\cite{OCAMPO2019}) are naturally realized in $VB_n/\Gamma_2(VP_n)$. 

Besides this, we shall prove that there are other elements of finite order (see Theorem~\ref{teo2}) and Bieberbach subgroups (see Theorem~\ref{teo6}) in $VB_n/\Gamma_2(VP_n)$.

\end{Obs}

The group $VB_n/\Gamma_2(VP_n)$ has a copy of the symmetric group $S_n$, since the short exact sequence~(\ref{seqvbn}) splits. 
In the next result, for $2\le t\le n$, we classify the elements of order $t$ in the crystallographic group $VB_n/\Gamma_2(VP_n)$ that projects onto a cycle element of $S_n$ of order $t$.


\begin{Teo}\label{teo2}
 
Consider $\sigma\in S_n$ of order $t$.
Let $\overline{\pi}_p(\sigma)^{-1}=\theta$ and let $T_{\theta}$ be a transversal of the action of $\sigma$ on $\{\lambda_{i,j}\mid 1\le i\neq j\le n\}$. 
The element $\prod_{1\le i\neq j\le n}\lambda_{i,j}^{a_{i,j}}\sigma$, where $a_{i,j}\in\Z$, has order $t$ in $VB_n/\Gamma_2(VP_n)$ if and only if $\sum_{\lambda_{i_1,j_1}\in \mathcal{O}_{\theta}(\lambda_{i,j})}a_{i_1,j_1}=0$ for all $\lambda_{i,j}\in T_{\theta}$.
\end{Teo}

\begin{proof}

Considering $\pi=\overline{\pi}_p$ and $\Z X=VP_n/\Gamma_2(VP_n)$ in Corollary~\ref{cor_teo2geral} we have the desired result. 

\end{proof}

Now we move to the conjugacy problem of elements (not necessarily of finite order) in the quotient group $VB_n/\Gamma_2(VP_n)$.
From Theorems~\ref{teo3} and~\ref{teo4} below, we will be able to determine whether any two elements in $VB_n/\Gamma_2(VP_n)$ are conjugate.

\begin{Teo}\label{teo3}

 Let $n\geq 2$ and let $\sigma\in VB_n/\Gamma_2(VP_n)$. Let $\theta=\overline{\pi}_p(\sigma)^{-1}$ and let $T_{\theta}$ be a transversal of the action of $\sigma$ on $\{\lambda_{i,j}\mid 1\le i\neq j\le n\}$. 
The elements $\prod_{1\le i\neq j\le n}\lambda_{i,j}^{a_{i,j}}\sigma$ and $\prod_{\lambda_{i_q,j_q}\in T_{\theta}}\lambda_{i_q,j_q}^{S_{i_q,j_q}}\sigma$ are conjugate in $VB_n/\Gamma_2(VP_n)$, where $S_{i_q,j_q}=\sum_{a_{r,s}\in\mathcal{O}_{\theta}(\lambda_{i_q,j_q})}{a_{r,s}}$ for all $\lambda_{i_q,j_q}\in T_{\theta}$.
\end{Teo}
\begin{proof}
In order to prove the theorem, just consider $G=VB_n/\Gamma_2(VP_n)$ and $\Z X=VP_n/\Gamma_2(VP_n)$ in Theorem~\ref{teo3_geral}.

\end{proof}

The next result states that the conjugacy classes of finite order elements in $VB_n/\Gamma_2(VP_n)$ are in correspondence with the conjugacy classes of elements in the symmetric group $S_n$, this is an analogous result for the case of the Artin braid group (see \cite[Theorem~5]{GGO1}).

\begin{Cor}~\label{conjordemfinita}
	Let $\beta_1,\beta_2$ elements of order $\tau$ in $VB_n/\Gamma_2(VP_n)$. Then $\beta_1$ and $\beta_2$ are conjugate in $VB_n/\Gamma_2(VP_n)$ if and only if $\overline{\pi}_p(\beta_1)$ and $\overline{\pi}_p(\beta_2)$ are conjugate in $S_n$.
\end{Cor}
\begin{proof}
It is clear that if $\beta_1$ and $\beta_2$ are conjugate in $VB_n/\Gamma_2(VP_n)$, then $\overline{\pi}_p(\beta_1)$ and $\overline{\pi}_p(\beta_2)$ are conjugate in $S_n$.

Now we prove the converse. Without loss of generality, we can assume that $\overline{\pi}_p(\beta_1)=\overline{\pi}_p(\beta_2)$. Let $\theta=\theta_1\theta_2\cdots\theta_r$ be a decomposition of $\overline{\pi}_p(\beta_1)$ and $\overline{\pi}_p(\beta_2)$ in $S_n$, where $\theta_i$ is a transposition. 
Conjugating $\overline{\pi}_p(\beta_1)$ and $\overline{\pi}_p(\beta_2)$ if necessary, we may suppose that $\overline{\pi}_p(\rho_{i_1}\rho_{i_2}\cdots\rho_{i_r})^{-1}=\overline{\pi}_p(\beta_1)=\overline{\pi}_p(\beta_2)$ where $\overline{\pi}_p(\rho_{i_j})^{-1}=\theta_j$. 
Thus, there are $X=\prod_{1\le i\neq j\le n}\lambda_{i,j}^{a_{i,j}}$ and $Y=\prod_{1\le i\neq j\le n}\lambda_{i,j}^{b_{i,j}}$ in $VP_n/\Gamma_2(VP_n)$ such that $\beta_1=X(\rho_{i_1}\rho_{i_2}\cdots\rho_{i_r})$ and $\beta_2=Y(\rho_{i_1}\rho_{i_2}\cdots\rho_{i_r})$. 
As $\beta_1$ and $\beta_2$ have finite order, we concluded by Theorem~\ref{teo2} that $\sum_{\lambda_{i_1,j_1}\in \mathcal{O}_{\theta}(\lambda_{i,j})}a_{i_1,j_1}=0$ and $\sum_{\lambda_{i_1,j_1}\in \mathcal{O}_{\theta}(\lambda_{i,j})}b_{i_1,j_1}=0$ for all $\lambda_{i,j}\in T_{\theta}$. 
Therefore, by Theorem~\ref{teo3}, $\beta_1$ and $\beta_2$ are conjugate in $VB_n/\Gamma_2(VP_n)$.
\end{proof}

In the following theorem, we study conjugacy classes of elements (not necessarily of finite order) in $VB_n/\Gamma_2(VP_n)$. 

\begin{Teo}\label{teo4}

Let $n\geq 2$ and let $\sigma\in VB_n/\Gamma_2(VP_n)$. Let $\theta=\overline{\pi}_p(\sigma)^{-1}$ and let $c\in C_{S_n}(\theta)$. 
Consider a transversal $T_{\theta}$ of the action of $\sigma$ on the set of generators $\{\lambda_{i,j}\mid 1\le i\neq j\le n\}$ of $VP_n/\Gamma_2(VP_n)$ and $\tilde{c}\in VB_n/\Gamma_2(VP_n)$ such that $\overline{\pi}_p(\tilde{c})^{-1}=c$ and $\tilde{c}\sigma\tilde{c}^{-1}=\prod_{1\le i\neq j}\lambda_{i,j}^{c_{i,j}}\sigma\in VP_n/\Gamma_2(VP_n)$. Then, the elements $\prod_{\lambda_{r,s}\in T_{\theta}}\lambda_{r,s}^{a_{r,s}}\sigma$ and $\prod_{\lambda_{r,s}\in T_{\theta}}\lambda_{r,s}^{b_{r,s}}\sigma$ are conjugate in $VB_n/\Gamma_2(VP_n)$ if and only if

\begin{enumerate}
\item[$(i)$] $\sum_{\lambda_{i,j}\in\mathcal{O}(\lambda_{r,s})}c_{i,j}+a_{r,s}=b_{r,s}$, for $\lambda_{c^{-1}(r),c^{-1}(s)}\in\mathcal{O}_{\theta}(\lambda_{r,s})$ where $\lambda_{r,s}\in T_{\theta}$.

\item[$(ii)$] $\sum_{\lambda_{i,j}\in\mathcal{O}(\lambda_{r,s})}c_{i,j}+a_{c^{-1}(r),c^{-1}(s)}=b_{r,s}$, for $\lambda_{c^{-1}(r),c^{-1}(s)}\not\in\mathcal{O}_{\theta}(\lambda_{r,s})$ where  $a_{c^{-1}(r),c^{-1}(s)}=\sum_{\lambda_{i,j}\in\mathcal{O}_{\theta}(\lambda_{c^{-1}(r),c^{-1}(s)})}a_{i,j}$ and $\lambda_{r,s}\in T_{\theta}$.
\end{enumerate}	
\end{Teo}

\begin{proof}
The proof of this result follows considering the free abelian group $\Z X=VP_n/\Gamma_2(VP_n)$ and the group $G=VB_n/\Gamma_2(VP_n)$ in Theorem~\ref{teo4_geral}. 
\end{proof}

Let $n\geq 2$. Let $\alpha\colon\Z\longrightarrow \aut(\Z_n)$ be a homomorphism, then $\alpha$ is determined by the image of a generator $x$ of $\Z$.
Further, there is $k\in\che 1,\dots,n-1\chd$ such that $\mdc(n,k)=1$ and  $\alpha(x)=\alpha_k$ where $\alpha_k\colon\Z_n\longrightarrow\Z_n$ is given by $\delta\longmapsto k\delta$.

\begin{Def}
A group $G$ is virtually cyclic if it has a cyclic subgroup $H$ of finite index.
\end{Def}
The above definition can be found in the introduction of \cite{FJ1993}.
In the next theorem, we realize virtually cyclic subgroups in $VB_n/\Gamma_2(VP_n)$.

\begin{Teo}\label{teo5}
 Let $n\ge 2$. Then, for every $k\in\{1,\dots,n-1\}$ such that $\mdc(n,k)=1$, the virtually cyclic group $\Z_n\rtimes_{\alpha_k}\Z$ can be realized in $VB_n/\Gamma_2(VP_n)$.
\end{Teo}
\begin{proof}
We show that there are elements $A$ and $B$ in $VB_n/\Gamma_2(VP_n)$ such that $A^n=1$, $BAB^{-1}=A^k$ and $\Z=\langle B\rangle$. Consider $A=\rho_1\rho_2\dots\rho_{n-1}$, thus, $A^n=1$. Suppose that there is $B\in VB_n/\Gamma_2(VP_n)$ such that $BAB^{-1}=A^k$, so, $\overline{\pi}_p(BAB^{-1})=\overline{\pi}_p(A^k)$. As by hypothesis $\mdc(n,k)=1$, the elements $\overline{\pi}_p(A)$ and $\overline{\pi}_p(A^k)$ are conjugate in $S_n$. Then there is $\gamma\in S_n$ so that $\gamma\overline{\pi}_p(A)\gamma^{-1}=\overline{\pi}_p(A^k)$. Consider $\tilde{\gamma}$ in $VB_n/\Gamma_2(VP_n)$ such that $\overline{\pi}_p(\tilde{\gamma})=\gamma$.

Let $B=\prod_{1\le i\neq j\le n}{\lambda_{i,j}^{x_{i,j}}}\tilde{\gamma}$. 
We will analyze which are the conditions on $x_{i,j}\in\Z$, $1\le i\neq j\le n$, so that $BAB^{-1}=A^k$ and $\Z=\langle B\rangle$. 
Since we want $\Z=\langle B\rangle$, then we consider $B$ in the group $VB_n/\Gamma_2(VP_n)$ of infinite order. On the other hand, $BAB^{-1}=A^k$ if and only if
\begin{align*}
BAB^{-1}=A^k&\Leftrightarrow\prod_{1\le i\neq j\le n}{\lambda_{i,j}^{x_{i,j}}}\tilde{\gamma}(\rho_1\rho_2\dots\rho_{n-1})\tilde{\gamma}^{-1}\prod_{1\le i\neq j\le n}{\lambda_{i,j}^{-x_{i,j}}}=(\rho_1\rho_2\dots\rho_{n-1})^k\\
&\Leftrightarrow\prod_{1\le i\neq j\le n}{\lambda_{i,j}^{x_{i,j}}}\tilde{\gamma}(\rho_1\rho_2\dots\rho_{n-1})\tilde{\gamma}^{-1}\prod_{1\le i\neq j\le n}{\lambda_{i,j}^{-x_{i,j}}}(\rho_1\rho_2\dots\rho_{n-1})^{-k}=1.
\end{align*}
Since $\overline{\pi}_p\pare\tilde{\gamma}(\rho_1\rho_2\dots\rho_{n-1})\tilde{\gamma}^{-1}\pard=(\rho_1\rho_2\dots\rho_{n-1})^k$, there are $c_{i,j}\in\Z$, $1\le i\neq j\le n$, such that $\tilde{\gamma}(\rho_1\rho_2\dots\rho_{n-1})\tilde{\gamma}^{-1}=\prod_{1\le i\neq j\le n}{\lambda_{i,j}^{c_{i,j}}}(\rho_1\rho_2\dots\rho_{n-1})^k$. So, considering $\theta=\overline{\pi}_p\pare(\rho_1\rho_2\dots\rho_{n-1})^k\pard^{-1}$ we have
\begin{align*}
	BAB^{-1}=A^k&\Leftrightarrow\prod_{1\le i\neq j\le n}{\lambda_{i,j}^{x_{i,j}}}\prod_{1\le i\neq j\le n}{\lambda_{i,j}^{c_{i,j}}}(\rho_1\rho_2\dots\rho_{n-1})^k\prod_{1\le i\neq j\le n}{\lambda_{i,j}^{-x_{i,j}}}(\rho_1\rho_2\dots\rho_{n-1})^{-k}=1\\
	&\Leftrightarrow\prod_{1\le i\neq j\le n}{\lambda_{i,j}^{x_{i,j}}}\prod_{1\le i\neq j\le n}{\lambda_{i,j}^{c_{i,j}}}\prod_{1\le i\neq j\le n}{\lambda_{\theta(i),\theta(j)}^{-x_{i,j}}}=1.
	\end{align*}
	Thus we obtain the following system of equations,
\begin{equation*}
	\begin{cases}
		x_{i,j}+c_{i,j}-x_{\theta^{-1}(i),\theta^{-1}(j)}=0,
	\end{cases}
\end{equation*}
 for all $1\le i\neq j\le n$. 
 Which is equivalent to the following subsystems of equations,
 \begin{equation}\label{sisciclicovirtual}
\begin{cases}
	x_{r,s}+c_{r,s}-x_{\theta^{-1}(r),\theta^{-1}(s)}=0\\
	x_{\theta^m(r),\theta^m(s)}+c_{\theta^m(r),\theta^m(s)}-x_{\theta^{m-1}(r),\theta^{m-1}(s)}=0
	\end{cases}
\end{equation}
for $\lambda_{r,s}\in T_{\theta}$ and $1\le m<|\mathcal{O}_{\theta}(\lambda_{r,s})|$. The system of equations~(\ref{sisciclicovirtual}) has solution if and only if, $c_{r,s}+c_{\theta(r),\theta(s)}+\cdots+c_{\theta^l(r),\theta^l(s)}=0$ for $l=|\mathcal{O}_{\theta}(\lambda_{r,s})|-1$ and $\lambda_{r,s}\in T_{\theta}$. 
Since $\prod_{1\le i\neq j\le n}{\lambda_{i,j}^{c_{i,j}}}(\rho_1\rho_2\dots\rho_{n-1})^k$ has finite order, therefore, by Theorem~\ref{teo2}, $c_{r,s}+c_{\theta(r),\theta(s)}+\cdots+c_{\theta^l(r),\theta^l(s)}=0$. So, the system of equations~(\ref{sisciclicovirtual}) has solution for all $\lambda_{r,s}\in T_{\theta}$. Thus, there is $B\in VB_n/\Gamma_2(VP_n)$ such that $BAB^{-1}=A^k$.
\end{proof}

\begin{Obs} 
The realization mentioned in Theorem~\ref{teo5} may be explicit. Taken the element $A$ given in the proof of that theorem, we may exhibit $\tilde{\gamma}\in VB_n/\Gamma_2(VP_n)$ as described in its proof. In order to do this, first write  $\tilde{\gamma}(\rho_1\dots\rho_{n-1})\tilde{\gamma}^{-1}(\rho_1\dots\rho_{n-1})^{-k}$ as a product of elements of $VP_n/\Gamma_2(VP_n)$ and then find one solution of the system of equations~(\ref{sisciclicovirtual}). After this procedure, we obtain an explicit realization of a virtually cyclic group $\Z_n\rtimes_{\alpha_k}\Z$ in $VB_n/\Gamma_2(VP_n)$.
\end{Obs}

\subsection{Realization of Bieberbach subgroups}\label{subsec:realizaBieberbach}

Torsion-free crystallographic groups, also called \textit{Bieberbach groups}, are notably interesting from its relation with Riemannian geometry and dynamical systems.  
It is well known (\cite[Section~3.3]{Wolf2011}) from the Bieberbach theorems that, via the fundamental group, there is a correspondence between the collection of Bieberbach groups of dimension $m$ and the collection of flat compact connected Riemannian manifolds of dimension $m$. 
A lot of information of a flat Riemannian manifold is recorded algebraically in its fundamental group (that is a Bieberbach group) or in its holonomy representation. 
So, it is very interesting to explore algebraic properties of Bieberbach groups in order to get information of their respective flat Riemannian manifolds. 
See for instance \cite{JR1991}, \cite{P1972},  \cite{22} and \cite{Wolf2011} for more details about these topics. 

We exhibit explicit Bieberbach subgroups of $VB_n/\Gamma_2(VP_n)$ with cyclic holonomy group in Theorem~\ref{teo6}. But first, we formulate the following lemma. 

\begin{Lema}\label{prodrho}
Let $n\geq 3$. The element $\rho_{r}\rho_{r+1}\cdots\rho_{r+t}$  has order ${t+2}$ in $VB_n$ for $1\le r\le n-2$, $t\ge 1$ and $r+t\le n-1$.
 \end{Lema}

\begin{proof} We recall that, in $S_n$, the element $\tau_{r}\tau_{r+1}\cdots\tau_{r+t}$ has order ${t+2}$ for $1\le r\le n-2$, $t\ge 1$ and $r+t\le n-1$.
The proof of this lemma follows from the fact that $VB_n=VP_n\rtimes S_n$ where the map $\iota\colon S_n\longrightarrow VB_n$ given by $\iota(\tau_i)=\rho_i$, for $i=1,\dots,n-1$ is a section for $\pi_p$.
\end{proof}

The proof of the next theorem is analogous to the one given for Theorem~3 of \cite{GGOP2021}. 

\begin{Teo}\label{teo6}
	Let $G_n$ be the subgroup of $S_n$ generated by $\pare 1,n,n-1,\dots,2\pard$. There is a Bieberbach subgroup  $\widetilde{G_n}$ in $VB_n/\Gamma_2(VP_n)$ of dimension $n(n-1)$ and whose holonomy group is $G_n$. Further, the center $Z(\widetilde{G_n})$ is a free Abelian group of rank $n-1$.
\end{Teo}

\begin{proof}
Let $\alpha_{n-1}=\rho_{n-1}\rho_{n-2}\dots\rho_{1}\in VB_n/\Gamma_2(VP_n)$, let $\tau=\overline{\pi}_p(\alpha_{n-1})^{-1}=\pare 1,n,n-1,\dots,2\pard$ and let $G_n=\langle\tau\rangle$. 
We know that $\alpha_{n-1}^n=1$ (see Lemma~\ref{prodrho}) in $VB_n/\Gamma_2(VP_n)$. 
Recall that the action of $\alpha_{n-1}$ on the basis $\{\lambda_{i,j}\mid 1\le i\neq j\le n\}$ of $VP_n/\Gamma_2(VP_n)$, induced from the action by conjugation of $VB_n/\Gamma_2(VP_n)$ on $VP_n/\Gamma_2(VP_n)$, is given by $\alpha_{n-1}\cdot \lambda_{i,j}=\lambda_{\tau(i), \tau(j)}$.
Then we have
\begin{align}\label{3.2}
(\lambda_{1,2}\alpha_{n-1})^n&=\prod_{\lambda_{r,s}\in\mathcal{O}_{\tau}(\lambda_{1,2})}\lambda_{r,s}.
\end{align} 
To simplify the notation, let $$
\prod_{\lambda_{r,s}\in\mathcal{O}_{\tau}(\lambda_{1,2})}\lambda_{r,s}=O_{1,2}.
$$ 
Consider the sets 
$$
X=\che\lambda_{1,2}\alpha_{n-1},\lambda_{i,j}^n\mid 1\le i\neq j\le n\chd,
$$
and 
$$
Y=\che O_{1,2},\lambda_{i,j}^n\mid 1\le i\neq j\le n\chd.
$$ 
Let $\widetilde{G_n}$ be the group generated by $X$ and let $L$ be the group generated by $Y$. Thus, the homomorphism ${\overline{\pi}_p}_{\arrowvert_{\widetilde{G_n}}}\colon\widetilde{G_n}\longrightarrow G_n$ is surjective and, by equation~(\ref{3.2}), $L$ is a subgroup of $\widetilde{G_n}$. We claim that $L=\ker \pare{\overline{\pi}_p}_{\arrowvert_{\widetilde{G_n}}}\pard$. Clearly, by the construction of $Y$, we have $L\subset\ker \pare{\overline{\pi}_p}_{\arrowvert_{\widetilde{G_n}}}\pard$. Conversely, let $w\in \ker \pare{\overline{\pi}_p}_{\arrowvert_{\widetilde{G_n}}}\pard$, we can write $w$ in the following way
$$
w=(\lambda_{1,2}\alpha_{n-1})^m\prod_{1\le i\neq j\le n}\lambda_{i,j}^{nt_{i,j}},
$$
where $m\in\Z$ and $t_{i,j}\in\Z$ for all $1\le i\neq j\le n$. 
Since $w\in \ker ({\overline{\pi}_p}_{\arrowvert_{\widetilde{G_n}}})$, it implies that $n$ divides $m$, so,
\begin{align*}
w&=\pare(\lambda_{1,2}\alpha_{n-1})^n\pard^{\frac{m}{n}}\prod_{1\le i\neq j\le n}\lambda_{i,j}^{nt_{i,j}}\\
&=O_{1,2}^{\frac{m}{n}}\prod_{1\le i\neq j\le n}\lambda_{i,j}^{nt_{i,j}}\in L.
\end{align*}
Thus, we have the following short exact sequence
\begin{equation*}
	\xymatrix{
		1 \ar[r] & L \ar[r]^-{} &  \widetilde{G_n} \ar[r]^-{{\overline{\pi}_p}_{\arrowvert_{\widetilde{G_n}}}} & G_n \ar[r] & 1.}
\end{equation*}
Notice that $L$ is a subgroup of $VP_n/\Gamma_2(VP_n)$ and 
$$
Y'=\che O_{1,2},\lambda_{i,j}^n\mid 1\le i\neq j\le n\mbox{ and }(i,j)\neq(1,2) \chd
$$ 
is a basis of the group $L$, since $\lambda_{1,2}^n=O_{1,2}^n\prod_{\lambda_{i,j}\in\mathcal{O}_{\tau}(\lambda_{1,2})}\lambda_{i,j}^{-n}$ where $(i,j)\neq (1,2)$. 
So, $L$ is a free Abelian group of rank $n(n-1)$. 
Using the action of $G_n$ on $L$, given by the restriction of the action of $S_n$ in $VP_n/\Gamma_2(VP_n)$, described in~(\ref{conjvbn}), we conclude that $\widetilde{G_n}$ is a crystallographic group of dimension $n(n-1)$ and whose holonomy is $G_n$. 
It remains to show that $\widetilde{G_n}$ is torsion free.

Let $w\in\widetilde{G_n}$ be an element of finite order. Since $L$ is torsion-free, the order of $w$ is equal to the order of $\overline{\pi}_p(w)$. In particular, $w^n=1$ since $\overline{\pi}_p(w)^n=1$. Using the action of $G_n$ on $L$, we guarantee the existence of $P\in L$ and $k\in\{1,2,\dots,n-1\}$ such that
$$
w=P(\lambda_{1,2}\alpha_{n-1})^k.
$$
We can write the element $P$ using the basis $Y'$, thereby we have 
$$
P=(O_{1,2})^{t_{1,2}}\prod_{1\le i\neq j\le n}\lambda_{i,j}^{nt_{i,j}},
$$
where $(i,j)\neq (1,2)$ and $t_{x,y}\in\Z$ for all $1\le x\neq y\le n$.
Thus,
\begin{align*}
1=w^n&=P(\lambda_{1,2}\alpha_{n-1})^kP(\lambda_{1,2}\alpha_{n-1})^k\cdots P(\lambda_{1,2}\alpha_{n-1})^k\\
&=P(\lambda_{1,2}\alpha_{n-1})^kP(\lambda_{1,2}\alpha_{n-1})^{-k}(\lambda_{1,2}\alpha_{n-1})^{2k}P(\lambda_{1,2}\alpha_{n-1})^{-2k}\\
&\cdots(\lambda_{1,2}\alpha_{n-1})^{(n-1)k}P(\lambda_{1,2}\alpha_{n-1})^{-(n-1)k}(\lambda_{1,2}\alpha_{n-1})^{nk}\\
&=P(\lambda_{1,2}\alpha_{n-1})^kP(\lambda_{1,2}\alpha_{n-1})^{-k}(\lambda_{1,2}\alpha_{n-1})^{2k}P(\lambda_{1,2}\alpha_{n-1})^{-2k}\\
&\cdots(\lambda_{1,2}\alpha_{n-1})^{(n-1)k}P(\lambda_{1,2}\alpha_{n-1})^{-(n-1)k}O_{1,2}^{k}.
\end{align*}
We notice that 
\begin{align*}
\pare\lambda_{1,2}\alpha_{n-1}\pard^kO_{1,2}\pare\lambda_{1,2}\alpha_{n-1}\pard^{-k}&=\pare\lambda_{1,2}\alpha_{n-1}\pard^k\prod_{\lambda_{r,s}\in\mathcal{O}_{\tau}(\lambda_{1,2})}\lambda_{r,s}\pare\lambda_{1,2}\alpha_{n-1}\pard^{-k}\\
&=\prod_{\lambda_{r,s}\in\mathcal{O}_{\tau}(\lambda_{1,2})}\lambda_{r,s}\\
&=O_{1,2}, 
\end{align*}
so, the exponent of $\lambda_{1,2}$ is equal to  $n\pare\sum_{\lambda_{r,s}\in\mathcal{O}_{\tau}(\lambda_{1,2})}t_{r,s}\pard+k$. 
We know that $w^n=1$ and $L$ is torsion-free, hence 
$$ 
n\pare\sum_{\lambda_{r,s}\in\mathcal{O}_{\tau}(\lambda_{1,2})}t_{r,s}\pard+k=0.
$$ 
Since $k\in\{1,2,\dots,n-1\}$ we have a contradiction. 
Therefore, does not exist $w$ in $\widetilde{G_n}$ of finite order. So $\widetilde{G_n}$ is a Bieberbach group of dimension $n(n-1)$ and with holonomy group $G_n$. 

To compute $Z(\widetilde{G_n})$ we will use Lemma~5.2 of \cite{22}, which ensures that the center of $\widetilde{G_n}$ is formed by the elements of $L$ fixed by the action of $G_n$. Using the basis $Y'$ of $L$ we conclude that
$$
\che\prod_{\lambda_{r,s}\in\mathcal{O}_{\tau}(\lambda_{1,2})}\lambda_{r,s},\prod_{\lambda_{x,y}\in\mathcal{O}_{\tau}(\lambda_{t,q})}\lambda_{x,y}^n\mid \lambda_{t,q}\in T_{\tau}\smallsetminus\mathcal{O}_{\tau}(\lambda_{1,2})\chd,
$$ 
is a basis for $Z(\widetilde{G_n})$. Note that the cardinality $\#T_{\tau}=n-1$ and therefore, $\Z(\widetilde{G_n})$ is free Abelian of rank $n-1$.
\end{proof}


Now, we recall the definitions of Anosov diffeomorphisms and K\"ahler structure of a given Riemannian manifold that will be used in the next theorem. 
First we note that in the cases of compact flat Riemannian manifolds these properties depend on completely of its holonomy representation, so we obtain geometric information about the manifold from the algebraic data of its fundamental group and holonomy representation.

A diffeomorphism $f\colon M \rightarrow M$ from a Riemannian manifold into itself	is called \textit{Anosov} and  $M$ is said to have a hyperbolic structure if the	following condition is satisfied: there exists a splitting of the tangent	bundle $T(M)=E^{s}+E^{u}$ such that $Df\colon E^{s} \rightarrow E^{s}$ is contracting and $Df\colon E^{u} \rightarrow E^{u}$ is expanding.
Classification of all compact manifolds supporting Anosov diffeomorphism play an important role in the theory of dynamical systems, that notion represents a kind of global hyperbolic behavior, and provides examples of stable dynamical systems. The problem of classifying those compact manifolds (up to diffeomorphism) which admit an Anosov diffeomorphism was first proposed by Smale \cite{SMALE1967}. 
An interesting algebraic characterization of the existence of Anosov diffeomorphisms, using the holonomy representation, was given by Porteous, see \cite[Theorems 6.1 and 7.1]{P1972}.	A  \textit{K\"ahler manifold} is a $2n$-real manifold endowed with Riemannian metric, complex structure, and a symplectic structure which are compatible at every point. 
An algebraic characterization of K\"ahler manifolds uses the knowledge of the behavior of certain irreducible representations that appear in the decomposition of the holonomy representation of the corresponding flat manifold, see \cite{JR1991} and \cite[Propositions~7.1~and~7.2]{22}.

\begin{Teo}\label{anosov}
Let $n\geq 2$, and let ${\mathcal M}_{n}$ be a 
$n(n-1)$-dimensional flat manifold whose fundamental group is the Bieberbach group $\widetilde{G_n}$ of Theorem~\ref{teo6}. Then, for every $n\geq 2$ ${\mathcal M}_{n}$ has first Betti number $n(n-1)$ and it is orientable if and only if $n$ is odd. In particular, for $n=2$ the 2-dimensional flat manifold ${\mathcal M}_{2}$ corresponds to the Klein bottle. 
For every $n\geq 3$ the flat manifold ${\mathcal M}_{n}$ admits Anosov diffeomorphisms, and it is a K\"ahler manifold if and only if $n$ is odd.

\end{Teo}

The proof of Theorem~\ref{anosov} is similar to the one given for \cite[Theorem~4]{GGOP2021} and it 
depends mainly on the holonomy representation of the Bieberbach group $\widetilde{G_n}$, analyzing the eigenvalues of the matrix representation and the decomposition of the holonomy representation in irreducible representations using character theory.

\begin{proof}

Let $n\geq 2$, and let ${\mathcal M}_{n}$ be a 
$n(n-1)$-dimensional flat manifold whose fundamental group is the Bieberbach group $\widetilde{G_n}$ of Theorem~\ref{teo6}. 
In terms of the basis $Y'=\che O_{1,2},\lambda_{i,j}^n\mid 1\le i\neq j\le n\mbox{ and }(i,j)\neq(1,2) \chd$ of the free Abelian group $L$, the holonomy representation  $\rho\colon G_n=\langle\tau\rangle \to \operatorname{\text{Aut}}(L)$ of the Bieberbach group $\widetilde{G_n}$ of Theorem~\ref{teo6} is given by the block diagonal matrix:
\begin{equation}\label{actmatrix}
\rho(\tau)=\left( \begin{smallmatrix}
M_{1} & & &\\
& M_{2} & &\\
& & \ddots &\\
& & & M_{n-1}
\end{smallmatrix}
\right),
\end{equation}
where $M_{1},\ldots,M_{n-1}$ are the $n$-by-$n$ matrices satisfying:
\begin{equation*}
\text{$M_{1}=\left(
\begin{smallmatrix}
1 & 0 & 0 & \cdots & 0 & 0 & n\\
0 & 0 & 0 & \cdots & 0 & 0 & -1\\
0 & 1 & 0 & \cdots & 0 & 0 & -1\\
0 & 0 & 1 & \cdots & 0 & 0 & -1\\
\vdots & \vdots & \vdots & \ddots & \vdots & \vdots & \vdots\\
0 & 0 & 0 & \cdots & 1 & 0 & -1\\
0 & 0 & 0 & \cdots & 0 & 1 & -1
\end{smallmatrix}\right)$ and $M_{2}=\cdots=M_{n-1}=
\left(
\begin{smallmatrix}
0 & 0 & 0 & \cdots & 0 & 0 & 1\\
1 & 0 & 0 & \cdots & 0 & 0 & 0\\
0 & 1 & 0 & \cdots & 0 & 0 & 0\\
0 & 0 & 1 & \cdots & 0 & 0 & 0\\
\vdots & \vdots & \vdots & \ddots & \vdots & \vdots & \vdots\\
0 & 0 & 0 & \cdots & 1 & 0 & 0\\
0 & 0 & 0 & \cdots & 0 & 1 & 0
\end{smallmatrix}
\right)$,}
\end{equation*}
where we have used the relation $\lambda_{1,2}^n=O_{1,2}^n\prod_{\lambda_{i,j}\in\mathcal{O}_{\tau}(\lambda_{1,2})}\lambda_{i,j}^{-n}$ with $(i,j)\neq (1,2)$.

By \cite[Theorems~6.4.12 and~6.4.13]{DEKIMPE1996}, the first Betti number of ${\mathcal M}_{n}$ is given by:
\begin{align*}
 \beta_{1}({\mathcal M}_{n})&=n(n-1) - \operatorname{\text{rank}}(\rho(\tau)- I_{n(n-1)})= n(n-1) - (n-1)(n-1)= n-1.
\end{align*}

Since, for every $i=1,\ldots, n-1$, the determinant of the matrix $M_i$ is equal to $(-1)^{n-1}$, then the determinant of the $n(n-1)$-by-$n(n-1)$ matrix $\rho(\tau)$ is 1 if and only if $n$ is odd. 
It is well known that a flat $m$-dimensional manifold determined by a Bieberbach group is orientable if and only if the image of its holonomy representation lies inside $SL(m, \Z)$ \cite[Theorem~6.4.6 and Remark~6.4.7]{DEKIMPE1996}. 
Then, for every $n\geq 2$ ${\mathcal M}_{n}$ is orientable if and only if $n$ is odd. 
Up to isomorphism, there are only two Bieberbach groups in dimension 2, one is the fundamental group of the torus and the other one corresponds to the fundamental group of the Klein bottle. 
Since $\widetilde{G_2}$ has holonomy $\Z_2$ then the 2-dimensional flat manifold ${\mathcal M}_{2}$ corresponds to the Klein bottle. 

Now we consider $n\geq 3$. We notice that, for every $i=1,\ldots, n-1$, the characteristic polynomial of $M_i$ is equal to $x^n-1$. Hence, by~(\ref{actmatrix}), the characteristic polynomial of $\rho(\tau)$ is equal to $(x^n-1)^{n-1}$, so the eigenvalues of $\rho(\tau)$ are the $n$th roots of unity each with multiplicity $n-1$. Therefore,  we conclude from \cite[Theorem~7.1]{P1972} that ${\mathcal M}_{n}$ admits Anosov diffeomorphisms. 

It remains to show when ${\mathcal M}_{n}$ admits a K\"ahler structure. In order to do this, we make use of the following result from~\cite[Theorem~3.1 and Proposition~3.2]{JR1991} that a Bieberbach group $\Pi$ of dimension $m$ is the fundamental group of a K\"ahler flat manifold with holonomy group $H$ if and only if $m$ is even, and each $\R$-irreducible summand of the holonomy representation $\psi\colon H \to \operatorname{\text{GL}}(m,\R)$ of $\Pi$, which is also $\mathbb{C}$-irreducible, occurs with even multiplicity. Since $\dim({\mathcal M}_{n})=n(n-1)$ and the character vector
of the representation $\rho$ is equal to $\left(\begin{smallmatrix}
n(n-1) \\
0\\
\vdots\\
0
\end{smallmatrix}\right)$,
it follows that each real irreducible representation of $\rho$ appears $n-1$ times in its decomposition, and hence that ${\mathcal M}_{n}$ admits a K\"ahler structure if and only if $n$ is odd.

\end{proof}

\subsection{The virtual braid group modulo the commutator subgroup of $KB_n$}\label{subsec:kbn}

Let $\pi_K\colon VB_n\longrightarrow S_n$ the homomorphism defined by $\pi_K(\sigma_i)=1$ and $\pi_K(\rho_i)=\tau_i$ for $i=1,\dots,n-1$ and $\tau_i=(i,i+1)$. The kernel of $\pi_K$ is denoted by $KB_n$, from which we obtain the following short exact sequence
\begin{equation*} 
	\xymatrix{
		1 \ar[r] & KB_n \ar[r]^-{} &  VB_n \ar[r]^-{{\pi_K}} &  S_n \ar[r] & 1}.
\end{equation*}

\begin{Prop}[Proposition~17 of \cite{BB2009}]
The group $KB_n$ admits a presentation with the generators $x_{k,l}$, $1\le k\neq l\le n$, and the defining relations
\begin{itemize}
    \item $x_{i,j}x_{k,l}=x_{k,l}x_{i,j}$,
    \item $x_{i,k}x_{k,j}x_{i,k}=x_{k,j}x_{i,k}x_{k,j}$,
\end{itemize}
where distinct letters stand for distinct indices. 
\end{Prop}

In this short subsection we deal with the quotient of the virtual braid group $VB_n$ modulo the commutator subgroup of $KB_n$, that fits into the following short exact sequence
\begin{equation}
	\xymatrix{
		1 \ar[r] & KB_n/\Gamma_2(KB_n) \ar[r]^-{} &  VB_n/\Gamma_2(KB_n) \ar[r]^-{\overline{\pi}_K} &   S_n \ar[r] & 1.}
\end{equation}

\begin{Prop}
The quotient group $VB_n/\Gamma_2(KB_n)$ is crystallographic if and only if $n=2$.
\end{Prop}

\begin{proof}

From Proposition~19 of~\cite{BB2009} we have that $KB_n/\Gamma_2(KB_n)$ is isomorphic to $\Z^2$ for $n=2, 3$ and it is isomorphic to $\Z$ for $n\ge 4$. 
The permutation action of $S_n$ on $KB_n$ was described in Proposition~3.1 of~\cite{BP2020} and it induces an action of $S_n$ on $KB_n/\Gamma_2(KB_n)$.

We note that $KB_2/\Gamma_2(KB_2)$ is generated by the cosets of $x_{1,2}$ and $x_{2,1}$ and clearly the holonomy representation $S_2\to \aut{(KB_2/\Gamma_2(KB_2))}$ is injective.

Let $n=3$. In the proof of Proposition~19 of~\cite{BB2009} the authors showed the equalities of cosets $x_{1,2}=x_{2,3}=x_{3,1}$ and $x_{1,3}=x_{3,2}=x_{2,1}$ in $KB_3/\Gamma_2(KB_3)$.
So, the automorphism of $KB_3/\Gamma_2(KB_3)=\Z^2$ induced by the permutation $(1,2,3)$ is the identity. 
Hence the representation $S_3\to \aut{(KB_3/\Gamma_2(KB_3))}$ is not injective and $VB_3/\Gamma_2(KB_3)$ is not crystallographic.

It is straightforward to see that the action of $S_n$ on $KB_n/\Gamma_2(KB_n)$ is not injective for $n\geq 4$, so we conclude that the group $VB_n/\Gamma_2(KB_n)$ is not a crystallographic.
\end{proof}

Although $VB_n/\Gamma_2(KB_n)$ is not crystallographic for $n\geq 3$, as we did in Theorems~\ref{teo2},~\ref{teo3} and~\ref{teo4}  for $VB_n/\Gamma_2(VP_n)$, applying the same techniques we may study structural aspects of the quotient $VB_n/\Gamma_2(KB_n)$. So we can state and prove equivalent versions of these theorems for $VB_n/\Gamma_2(KB_n)$. For instance, we illustrate it for the particular case of 3 strings in the followings result.

\begin{Teo}

\begin{enumerate}
    \item Let $\rho_{r_1}\in VB_3/\Gamma_2(KB_3)$. The element $x_{1,2}^{a_{1,2}}x_{1,3}^{a_{1,3}}\rho_{r_1}$, where $a_{i,j}\in\Z$, have order $2$ in $VB_3/\Gamma_2(KB_3)$ if and only if $a_{1,2}+a_{1,3}=0$.

    \item Let $\rho_{r_1}\in VB_3/\Gamma_2(KB_3)$. The elements $x_{1,2}^{a_{1,2}}x_{1,3}^{a_{1,3}}\rho_{r_1}$ and $x_{1,2}^{b_{1,2}}x_{1,3}^{b_{1,3}}\rho_{r_1}$ are conjugate in $VB_3/\Gamma_2(KB_3)$ if and only if $a_{1,2}+a_{1,3}=b_{1,2}+b_{1,3}$, where $a_{i,j},b_{i,j}\in\Z$.

\end{enumerate}
\end{Teo}

\section{Crystallographic groups and virtual twin groups}\label{groupvtn}

The connection between virtual twin groups and crystallographic groups will be treated in this section. 
The results stated in this section are analogous versions of the ones given in Section~\ref{groupvbn} for the quotient group $VT_n/\Gamma_2(PVT_n)$. 
We start rewriting a presentation of the virtual twin group $VT_n$, it was formulated in \cite[p.3]{NNS2020}.

\begin{Teo}[\cite{NNS2020}]\label{apvtn}
The {virtual twin group} $VT_n$, $n\ge 2$, admits the following group presentation
\begin{itemize}
    \item Generators:
    \subitem $\sigma_i,\rho_i$ for $i=1,\dots,n-1$.
\item Relations:
\subitem $\sigma_i^2=1$ for $i=1,2,\dots,n-1$.
\subitem $\sigma_i\sigma_j=\sigma_j\sigma_i$ for $|i-j|\ge 2$.
\subitem $\rho_i^2=1$ for $i=1,\dots,n-1$.
\subitem $\rho_i\rho_j=\rho_j\rho_i$ for $|i-j|\ge 2$.
\subitem $\rho_i\rho_{i+1}\rho_{i}=\rho_{i+1}\rho_{i}\rho_{i+1}$, for $i=1,2,\dots,n-2$.
\subitem $\rho_i\sigma_j=\sigma_j\rho_i$, for $|i-j|\ge 2$.
\subitem $\rho_i\rho_{i+1}\sigma_{i}=\sigma_{i+1}\rho_{i}\rho_{i+1}$, for $i=1,\dots,n-2$.
\end{itemize}
\end{Teo}

 Let $n\geq 2$ and let $\pi\colon VT_n\longrightarrow S_n$ be the homomorphism defined by ${\pi}(\sigma_i)={\pi}(\rho_i)=\tau_i$ for $i=1,\dots,n-1$ and $\tau_i=(i,i+1)$. The \textit{pure virtual twin group} $PVT_n$ is defined to be the kernel of $\pi$, from which we obtain the following short exact sequence
\begin{equation}\label{seqvtn}  
	\xymatrix{
		1 \ar[r] & PVT_n \ar[r]^-{} &  VT_n \ar[r]^-{{\pi}} &  S_n \ar[r] & 1.}
\end{equation}

As mentioned at the top of page 5 of \cite{NNS2020}, the  virtual twin group admits a decomposition as semi product $VT_n=PVT_n\rtimes S_n$, the map $\iota\colon S_n\longrightarrow VT_n$ given by $\iota(\tau_i)=\rho_i$, for $i=1,\dots,n-1$, is naturally a section for $\pi$.
%
%
We state here, without proof, the following presentation of the pure virtual twin group.

\begin{Teo}[Theorem~3.3, \cite{NNS2020}]\label{appvtn}
The pure virtual twin group $PVT_n$ on $n\ge 2$ strands has the presentation
\begin{itemize}
    \item Generators:
    \subitem $\lambda_{i,j}$ for $1\le i<j\le n$.
    \item Relations:
    \subitem $\lambda_{i,j}\lambda_{k,l}=\lambda_{k,l}\lambda_{i,j}$ for distinct integers $i,j,k,l$.
\end{itemize}
\end{Teo}

\begin{Obs}
In the presentation of $PVT_n$ given in \cite{NNS2020} the element $\lambda_{i,i+1}$ is equal to $\sigma_i\rho_i$, for $i=1,\dots,n-2$. 
On the other hand, in the presentation of $VP_n$ given in \cite{B2004} the element  $\lambda_{i,i+1}$ is equal to $\rho_i\sigma_i^{-1}$, for $i=1,\dots,n-2$. 
We note that  $(\sigma_i\rho_i)^{-1}=\rho_i\sigma_i^{-1}$ for $i=1,\dots,n-2$.
\end{Obs}

From Remark~3.2 of~\cite{NNS2020}, using the presentation of the pure virtual twin group given in Theorem~\ref{appvtn}, we conclude that the action by conjugation of $S_n=\langle\rho_1,\dots,\rho_{n-1}\rangle$ on
$PVT_n$ is given by
\begin{equation}\label{actionvtn}
 \rho\cdot\lambda_{i,j}=\lambda_{\rho^{-1}(i),\rho^{-1}(j)},  
\end{equation}
where $\lambda_{j,i}=\lambda_{i,j}^{-1}$ for all $1\le i<j\le n$.

Recall that the pure virtual twin group $PVT_n$ is defined to be the kernel of $\pi\colon VT_n \to S_n$. The split short exact sequence (\ref{seqvtn}) provides, for every $n\geq 2$, the following short exact sequence
\begin{equation}\label{seqvtn/} 
	\xymatrix{
		1 \ar[r] & PVT_n/\Gamma_2(PVT_n) \ar[r]^-{} &  VT_n/\Gamma_2(PVT_n) \ar[r]^-{\overline{\pi}} &   S_n \ar[r] & 1.}
\end{equation}
We notice that $\iota$, the section for $\pi$ described in the paragraph below equation~(\ref{seqvtn}), induces a homomorphism $\overline{\iota}\colon S_n\longrightarrow VT_n/\Gamma_2(PVT_n)$ which is a section for $\overline{\pi}$, so we obtain the decomposition $VT_n/\Gamma_2(PVT_n)=PVT_n/\Gamma_2(PVT_n)\rtimes S_n$.

Using the same techniques applied in Section~\ref{groupvbn} we may obtain several results about the structure of $VT_n/\Gamma_2(PVT_n)$. 
For instance, we will show that the group $VT_n/\Gamma_2(PVT_n)$ is a crystallographic group, we will characterize the elements of finite order and the conjugacy classes of elements in $VT_n/\Gamma_2(PVT_n)$, and also we will realize Bieberbach groups and some infinite virtually cyclic groups in $VT_n/\Gamma_2(PVT_n)$. 

The first result in this section is that $VT_n/\Gamma_2(PVT_n)$ is a crystallographic group of dimension $n(n-1)/2$ and whose holonomy group is $S_n$.

\begin{Teo}\label{teo7}
Let $n\ge 2$.  There is a split short exact sequence:
	\begin{equation}\label{seqvtn2} 
	\xymatrix{
		1 \ar[r] & \Z^{n(n-1)/2} \ar[r]^-{} &  VT_n/\Gamma_2(PVT_n) \ar[r]^-{\overline{\pi}} &   S_n \ar[r] & 1}
		\end{equation}
	and the middle group is a crystallographic group.
\end{Teo}

\begin{proof}
From the presentation of the pure virtual twin group given in Theorem~\ref{appvtn} it is straightforward to show that $PVT_n/\Gamma_2(PVT_n)$ is isomorphic to the free Abelian group $\Z^{n(n-1)/2}$ generated by the set 
\begin{equation}\label{basispvt}
\{\lambda_{i,j}\mid 1\le i< j\le n\}. 
\end{equation}
The action of $S_n$ on $PVT_n/\Gamma_2(PVT_n)$ is induced from the action of $S_n$ on $PVT_n$ (see equation~(\ref{actionvtn})). 
So the proof follows in an analogous way to the proof of  Theorem~\ref{teo1}. 

\end{proof}

The group $ VT_n/\Gamma_2(PVT_n)$ has a copy of the symmetric group $S_n$, since the short exact sequence~(\ref{seqvtn/}) splits. It is natural to ask which other elements in this group have finite order. We classify the elements of order $t$, for $2\le t\le n$, in the crystallographic group $VT_n/\Gamma_2(PVT_n)$ that project to a cycle element of $S_n$ of order $t$. 

\begin{Teo}\label{teo8}
Consider $\sigma\in S_n$ of order $t$. Let $\overline{\pi}(\sigma)^{-1}=\theta$ and let $T_{\theta}$ be a transversal of the action of $\sigma$ on $\{\lambda_{i,j}\mid 1\le i< j\le n\}$. The element $\prod_{1\le i< j\le n}\lambda_{i,j}^{a_{i,j}}\sigma$, where $a_{i,j}\in\Z$, have order $t$ in $VT_n/\Gamma_2(PVT_n)$ if and only if for all $\lambda_{r,s}\in T_{\theta}$ 
we have
\begin{align*}
a_{r,s}+(-1)^{k_{\theta^{-1}(r),\theta^{-1}(s)}}a_{\theta^{-1}(r),\theta^{-1}(s)}+(-1)^{k_{\theta^{-2}(r),\theta^{-2}(s)}}a_{\theta^{-2}(r),\theta^{-2}(s)}+\dots +\\
(-1)^{k_{\theta^{-q}(r),\theta^{-q}(s)}}a_{\theta^{-q}(r),\theta^{-q}(s)}=0,
\end{align*}
where $\lambda_{\theta^i(r),\theta^i(s)}\in\mathcal{O}_{\theta}(\lambda_{r,s})$, $k_{\theta^i(r),\theta^i(s)}=0$ if $1\le \theta^i(r)<\theta^i(s)\le n$ and $k_{\theta^i(r),\theta^i(s)}=1$, $a_{\theta^i(r),\theta^i(s)}=a_{\theta^i(s),\theta^i(r)}$  if $1\le \theta^i(s)<\theta^i(r)\le n$.
\end{Teo}

\begin{proof}
Taking $G=VT_n/\Gamma_2(PVT_n)$ and $\Z X=PVT_n/\Gamma_2(PVT_n)$  in Corollary~\ref{cor_teo2geral} we get the proof.

\end{proof}


Now we consider the conjugacy problem of elements, not necessarily of finite order, in the quotient group $VT_n/\Gamma_2(PVT_n)$. 

\begin{Teo}\label{teo9}
Let $\sigma\in VT_n/\Gamma_2(PVT_n)$, let $\theta=\overline{\pi}(\sigma)^{-1}$ and let $T_{\theta}$ be a transversal of the action of $\sigma$ on $\{\lambda_{i,j}\mid 1\le i< j\le n\}$. The elements $\prod_{1\le i< j\le n}\lambda_{i,j}^{a_{i,j}}\sigma$ and $\prod_{\lambda_{r,s}\in T_{\theta}}\lambda_{r,s}^{b_{r,s}}\sigma$ are conjugate in $VT_n/\Gamma_2(PVT_n)$, if  for all $\lambda_{r,s}\in T_{\theta}$ we have
$$
b_{r,s}=
\begin{cases}
a_{r,s}+(-1)^{k_{r_1,s_1}}a_{r_1,s_1}+\dots +(-1)^{k_{r_k,s_k}}a_{r_k,s_k},\mbox{ if }\\ \mathcal{O}_{\theta}(\lambda_{r,s})\neq\che\lambda_{r,s},\lambda_{r,s}^{-1},\dots,\lambda_{r_l,s_l},\lambda_{r_l,s_l}^{-1}\chd,\\

a_{r,s}-a_{\theta(r),\theta(s)}-\dots-a_{\theta^q(r),\theta^q(s)},\mbox{ if } \mathcal{O}_{\theta}(\lambda_{r,s})= \che\lambda_{r,s},\lambda_{r,s}^{-1},\dots,\lambda_{r_l,s_l},\lambda_{r_l,s_l}^{-1}\chd
\end{cases}
$$
where $q=\frac{|\mathcal{O}_{\theta}(\lambda_{r,s})|}{2}-1$, $k_{\theta^i(r),\theta^i(s)}=0$ if $1\le \theta^i(r)<\theta^i(s)\le n$, and $k_{\theta^i(r),\theta^i(s)}=1$ with $a_{\theta^i(r),\theta^i(s)}=a_{\theta^i(s),\theta^i(r)}$  if $1\le \theta^i(s)<\theta^i(r)\le n$.
\end{Teo}

\begin{proof}
To prove the result, consider $G=VT_n/\Gamma_2(PVT_n)$ and $\Z X=PVT_n/\Gamma_2(PVT_n)$ in Theorem~\ref{teo3_geral}.

\end{proof}

Just as we have done in Corollary~\ref{conjordemfinita}, we can show that the conjugacy classes of finite order elements in $VT_n/\Gamma_2(PVT_n)$ are in correspondence with the conjugacy classes of elements in the symmetric group $S_n$. 
Theorem~\ref{teo10} is the analogous version  of Theorem~\ref{teo4} but now for the quotient of the virtual twin group $VT_n/\Gamma_2(PVT_n)$. 
We recall that $C_{S_n}(\theta)$ denotes the centralizer of $\theta$ in $S_n$.

\begin{Teo}\label{teo10}

Let $\sigma\in VT_n/\Gamma_2(PVT_n)$, let $\theta=\overline{\pi}(\sigma)^{-1}$, let $c\in C_{S_n}(\theta)$ and let $T_{\theta}$ be a transversal of the action of $\sigma$ on $\{\lambda_{i,j}\mid 1\le i< j\le n\}$. 
Consider $\tilde{c}\in VT_n/\Gamma_2(PVT_n)$ such that $\overline{\pi}(\tilde{c})^{-1}=c$ and $\tilde{c}\sigma\tilde{c}^{-1}=\prod_{1\le i< j\le n}\lambda_{i,j}^{c_{i,j}}\sigma\in VT_n/\Gamma_2(PVT_n)$. The elements $\prod_{\lambda_{r,s}\in T_{\theta}}\lambda_{r,s}^{a_{r,s}}\sigma$ and $\prod_{\lambda_{r,s}\in T_{\theta}}\lambda_{r,s}^{b_{r,s}}\sigma$ are conjugate in $VT_n/\Gamma_2(PVT_n)$ if and only if
\begin{itemize}
\item[(1)] For all $\lambda_{r,s}\in T_{\theta}$ such that $\mathcal{O}_{\theta}(\lambda_{r,s})\neq\che\lambda_{r,s},\lambda_{r,s}^{-1},\dots,\lambda_{r_k,s_k},\lambda_{r_k,s_k}^{-1}\chd$ we have
\begin{itemize}
\item[(i)] $\sum_{\lambda_{i,j}\in\mathcal{O}(\lambda_{r,s})}c_{i,j}+(-1)^{k_{r,s}}a_{r,s}=b_{r,s}$, if $\lambda_{c^{-1}(r),c^{-1}(s)}\in\mathcal{O}_{\theta}(\lambda_{r,s})$,
    
\item[(ii)] $\sum_{\lambda_{i,j}\in\mathcal{O}(\lambda_{r,s})}c_{i,j}+(-1)^{k_{r,s}}a_{c^{-1}(r),c^{-1}(s)}=b_{r,s}$, if $\lambda_{c^{-1}(r),c^{-1}(s)}\notin\mathcal{O}_{\theta}(\lambda_{r,s})$ for  $\lambda_{r,s}\in T_{\theta}$ and
$$
a_{c^{-1}(r),c^{-1}(s)}=\sum_{\lambda_{i,j}\in\mathcal{O}_{\theta}(\lambda_{c^{-1}(r),c^{-1}(s)})}a_{i,j},
$$ 

\end{itemize}
 where $k_{r,s}=0$ if $c(r)<c(s)$ and $k_{r,s}=1$ if $c(r)>c(s)$.
\item[(2)] For all $\lambda_{r,s}\in T_{\theta}$ such that $\mathcal{O}_{\theta}(\lambda_{r,s})=\che\lambda_{r,s},\lambda_{r,s}^{-1},\dots,\lambda_{r_k,s_k},\lambda_{r_k,s_k}^{-1}\chd$ we have $2k+a_{r,s}=b_{r,s}$ for $k\in\Z$.
\end{itemize}
\end{Teo}

\begin{proof}
In order to prove the result, consider in Theorem~\ref{teo4_geral} the groups $G=VT_n/\Gamma_2(PVT_n)$ and $\Z X=PVT_n/\Gamma_2(PVT_n)$.

\end{proof}

Let $n\geq 2$. Let $\alpha\colon\Z\longrightarrow \aut(\Z_n)$ be a homomorphism, then $\alpha$ is determined by the image of a generator $x$ of $\Z$. Recall that, there is $k\in\che 1,\dots,n-1\chd$ such that $\mdc(n,k)=1$ and  $\alpha(x)=\alpha_k$ where $\alpha_k\colon\Z_n\longrightarrow\Z_n$ is given by $\delta\longmapsto k\delta$. 

\begin{Teo}\label{teo11}
 Let $n\ge 2$. Then, for every $k\in\{1,\dots,n-1\}$ such that $\mdc(n,k)=1$, the virtually cyclic group $\Z_n\rtimes_{\alpha_k}\Z$ can be realized in $VT_n/\Gamma_2(PVT_n)$.
\end{Teo}

\begin{proof}
The proof of this theorem is analogous to the one given for Theorem~\ref{teo5}, but using the action of $S_n$ on $PVT_n$ described in equation~(\ref{actionvtn}) and Theorem~\ref{teo8}.
\end{proof}

Finally, we note that  following the idea of the proof of Theorem~\ref{teo6}, we can exhibit explicit Bieberbach subgroups, with cyclic holonomy group, in the quotient of the virtual twin group $VT_n/\Gamma_2(PVT_n)$. 
So, there exists compact flat Riemannian manifolds such that its fundamental group is a subgroup of  $VT_n/\Gamma_2(PVT_n)$ and for those manifolds we may obtain several geometric information using the algebraic data in its fundamental group (see Theorem~\ref{anosov} and its proof).

\section{On quotients of  virtual braid groups and virtual twin groups}\label{quotgroupvbn}

In the first subsection we consider some quotients of virtual braid groups and virtual twin groups and discuss if we may relate them to crystallographic groups as we did with $VB_n$ and $VT_n$ in previous sections.
In the second subsection, we deal with the extended pure loop braid group. 


\subsection{Welded, unrestricted, flat virtual and flat welded braid groups}\label{subsec:wf}



In this subsection we shall study many of the quotients of the virtual braid group $VB_n$ and of the virtual twin group given in Figure~\ref{quocientes}. 
For the sake of completeness, we recall here the corresponding definitions.

The \textit{welded braid group} $WB_n$ can be defined (see \cite[p.4]{BBD2015}) as a quotient of $VB_n$ by the normal subgroup generated by the relations
$$
\rho_i\sigma_{i+1}\sigma_i=\sigma_{i+1}\sigma_{i}\rho_{i+1},\mbox{ for }i=1,\dots,n-2.
$$
The \textit{unrestricted virtual braid} $UVB_n$ can be defined (see Definition~2.3 of \cite{BBD2015}) as a quotient of $WB_n$ 
by the normal subgroup generated by the relations
$$
\rho_{i+1}\sigma_{i}\sigma_{i+1}=\sigma_{i}\sigma_{i+1}\rho_{i}\mbox{ for }i=1,2,\dots,n-2.
$$


Let $\pi_P\colon WB_n\longrightarrow S_n$
 (resp. $\pi_P\colon UVB_n\longrightarrow S_n$), where ${\pi_P}(\sigma_i)={\pi}(\rho_i)=\tau_i$ for $i=1,\dots,n-1$ and $\tau_i=(i,i+1)$. The \textit{pure welded braid group} $WP_n$ (resp. \textit{unrestricted virtual pure braid group} $UVP_n$) is defined to be the kernel of $\pi_P$, from which we obtain the following short exact sequence
 \begin{equation*} 
	\xymatrix{
		1 \ar[r] & WP_n \ar[r]^-{} &  WB_n \ar[r]^-{{\pi}_p} &  S_n \ar[r] & 1,}
\end{equation*}
  
\begin{equation*} 
	(\textrm{resp. }\,
	\xymatrix{
		1 \ar[r] & UVP_n \ar[r]^-{} &  UVB_n \ar[r]^-{{\pi_P}} &  S_n \ar[r] & 1,})
\end{equation*}
that provides us the following short exact sequence
\begin{equation}
	\xymatrix{
		1 \ar[r] & WP_n/\Gamma_2(WP_n) \ar[r]^-{} &  WB_n/\Gamma_2(WP_n) \ar[r]^-{\overline{\pi}_p} &   S_n \ar[r] & 1.}
\end{equation}
\begin{equation}
	(\textrm{resp. }\,
	\xymatrix{
		1 \ar[r] & UVP_n/\Gamma_2(UVP_n) \ar[r]^-{} &  UVB_n/\Gamma_2(UVP_n) \ar[r]^-{\overline{\pi_P}} &   S_n \ar[r] & 1.})
\end{equation}
Presentations for these groups may be found in Corollary~3.19 of~\cite{D2017} for $WP_n$ and in Theorem~2.7 of \cite{BBD2015} for $UVP_n$.

In the following result we show that the quotient groups $WB_n/\Gamma_2(WP_n)$, $UVB_n/\Gamma_2(UVP_n)$ and $VB_n/\Gamma_2(VP_n)$ are isomorphic. 
So, it is a consequence that all of them are crystallographic groups, and the results given in Section~\ref{groupvbn} also holds for $WB_n/\Gamma_2(WP_n)$ and $UVB_n/\Gamma_2(UVP_n)$.

\begin{Teo}\label{teo:isowbn/gamma}
Let $n\geq 2$. 
The groups $WB_n/\Gamma_2(WP_n)$ and $UVB_n/\Gamma_2(UVP_n)$ are isomorphic to $VB_n/\Gamma_2(VP_n)$.
\end{Teo}

\begin{proof}
We need to show that $\sigma_{i+1}\sigma_{i}\rho_{i+1}\sigma_i^{-1}\sigma_{i+1}^{-1}\rho_{i}^{-1}$ and $\rho_{i+1}\sigma_{i}\sigma_{i+1}\rho_{i}^{-1}\sigma_{i+1}^{-1}\sigma_{i}^{-1}$ are elements of $\Gamma_2(VP_n)$ for all $i=1,2,\dots,n-2$ (see Figure~\ref{quocientes}). 
We only prove that $\sigma_{i+1}\sigma_{i}\rho_{i+1}\sigma_i^{-1}\sigma_{i+1}^{-1}\rho_{i}^{-1}\in\Gamma_2(VP_n)$, since the proof of the claim $\rho_{i+1}\sigma_{i}\sigma_{i+1}\rho_{i}^{-1}\sigma_{i+1}^{-1}\sigma_{i}^{-1}\in\Gamma_2(VP_n)$ can be verified similarly. 
We use the relations given in Theorem~\ref{apvbn} to verify it:
\begin{align*}
\sigma_{i+1}\sigma_i\rho_{i+1}\sigma_{i}^{-1}\sigma_{i+1}^{-1}\rho_i^{-1}&=\rho_{i+1}\rho_{i+1}\sigma_{i+1}\sigma_i\rho_{i+1}\sigma_{i}^{-1}\sigma_{i+1}^{-1}\rho_i^{-1}\\    
&=\rho_{i+1}\lambda_{i+2,i+1}^{-1}\sigma_i\rho_{i+1}\sigma_i^{-1}\sigma_{i+1}^{-1}\rho_i^{-1}\\
&=\rho_{i+1}\lambda_{i+2,i+1}^{-1}\rho_i\rho_i\sigma_i\rho_{i+1}\sigma_{i}^{-1}\sigma_{i+1}^{-1}\rho_i^{-1}\\
&=\rho_{i+1}\lambda_{i+2,i+1}^{-1}\rho_i\lambda_{i+1,i}^{-1}\rho_{i+1}\sigma_i^{-1}\sigma_{i+1}^{-1}\rho_i^{-1}\\
&=\rho_{i+1}\lambda_{i+2,i+1}^{-1}\rho_i\lambda_{i+1,i}^{-1}\rho_{i+1}\sigma_i^{-1}\rho_i\rho_i\sigma_{i+1}^{-1}\rho_i^{-1}\\
&=\rho_{i+1}\lambda_{i+2,i+1}^{-1}\rho_i\lambda_{i+1,i}^{-1}\rho_{i+1}\lambda_{i+1,i}\rho_i\sigma_{i+1}^{-1}\rho_i^{-1}\\
&=\rho_{i+1}\lambda_{i+2,i+1}^{-1}\rho_i\lambda_{i+1,i}^{-1}\rho_{i+1}\lambda_{i+1,i}\rho_i\lambda_{i+2,i+1}\rho_{i+1}\rho_i^{-1}\\
&=\rho_{i+1}\lambda_{i+2,i+1}^{-1}\lambda_{i,i+1}^{-1}\rho_i\rho_{i+1}\rho_i\lambda_{i,i+1}\lambda_{i+2,i+1}\rho_{i+1}\rho_i\\
&=\rho_{i+1}\lambda_{i+2,i+1}^{-1}\lambda_{i,i+1}^{-1}\rho_{i+1}\rho_i\rho_{i+1}\lambda_{i,i+1}\lambda_{i+2,i+1}\rho_{i+1}\rho_i\\
&=\lambda_{i+1,i+2}^{-1}\lambda_{i,i+2}^{-1}\rho_i\lambda_{i,i+2}\lambda_{i+1,i+2}\rho_i\\
&=\lambda_{i+1,i+2}^{-1}\lambda_{i,i+2}^{-1}\lambda_{i+1,i+2}\lambda_{i,i+2}.
\end{align*}

Thereby, $\sigma_{i+1}\sigma_{i}\rho_{i+1}\sigma_i^{-1}\sigma_{i+1}^{-1}\rho_{i}^{-1}\in\Gamma_2(VP_n)$ for all $i=1,2,\dots,n-2$.
\end{proof}



The \textit{flat virtual braid group} $FVB_n$ can be defined (see \cite[p.17]{BBD2015}) as a quotient of $VB_n$. 
However, using the presentations of the groups (see Figure~\ref{quocientes}), we note that it is straightforward to show that $FVB_n$ can also be viewed as a quotient of 
$VT_n$ by the normal subgroup generated by the relations
$$
\sigma_i\sigma_{i+1}\sigma_i=\sigma_{i+1}\sigma_{i}\sigma_{i+1},\mbox{ for }i=1,\dots,n-2.
$$
Similarly, the \textit{flat welded braid group} $FWB_n$ can be defined as a quotient of $WB_n$ (see Section~5.2 of~\cite{BBD2015}). 
Also, as described in Figure~\ref{quocientes}, $FWB_n$ is the quotient of $FVB_n$ by the normal subgroup generated by the relations  
$$
\rho_i\sigma_{i+1}\sigma_i=\sigma_{i+1}\sigma_{i}\rho_{i+1},\mbox{ for }i=1,\dots,n-2
$$
or, as the quotient of $WT_n$ by the normal subgroup generated by the relations
$$
\sigma_i\sigma_{i+1}\sigma_i=\sigma_{i+1}\sigma_{i}\sigma_{i+1},\mbox{ for }i=1,\dots,n-2.
$$

Let $\pi_p\colon FWB_n\longrightarrow S_n$
 (resp. $\pi_p\colon FVB_n\longrightarrow S_n$), where ${\pi_p}(\sigma_i)={\pi_p}(\rho_i)=\tau_i$ for $i=1,\dots,n-1$ and $\tau_i=(i,i+1)$. The \textit{flat pure welded braid group} $FWP_n$ (resp. \textit{flat virtual pure braid group} $FVP_n$) is defined to be the kernel of $\pi_p$, from which we obtain the following short exact sequence
 \begin{equation*} 
	\xymatrix{
		1 \ar[r] & FWP_n \ar[r]^-{} &  FWB_n \ar[r]^-{{\pi}_p} &  S_n \ar[r] & 1,}
\end{equation*}
  
\begin{equation*} 
	(\textrm{resp. }\,
	\xymatrix{
		1 \ar[r] & FVP_n \ar[r]^-{} &  FVB_n \ar[r]^-{{\pi}_p} &  S_n \ar[r] & 1,})
\end{equation*}
that provides us the following short exact sequence
\begin{equation}
	\xymatrix{
		1 \ar[r] & FWP_n/\Gamma_2(FWP_n) \ar[r]^-{} &  FWB_n/\Gamma_2(FWP_n) \ar[r]^-{\overline{\pi}_p} &   S_n \ar[r] & 1.}
\end{equation}
\begin{equation}
	(\textrm{resp. }\,
	\xymatrix{
		1 \ar[r] & FVP_n/\Gamma_2(FVP_n) \ar[r]^-{} &  FVB_n/\Gamma_2(FVP_n) \ar[r]^-{\overline{\pi}_p} &   S_n \ar[r] & 1.})
\end{equation}
Presentations for these groups may be found in  the paragraph before Proposition~5.5 of~\cite{BBD2015} for $FWP_n$ and in Proposition~5.1 of \cite{BBD2015} for $FVP_n$.

The next statement is similar to Theorem~\ref{teo:isowbn/gamma}, but now related to virtual twin groups. 
We show that the quotient groups $FWB_n$, $FVB_n/\Gamma_2(FVP_n)$ and $VT_n/\Gamma_2(PVT_n)$ are isomorphic. 
Hence, all of them are crystallographic groups and all results given in Section~\ref{groupvtn} also holds for $FWB_n$ and $FVB_n/\Gamma_2(FVP_n)$.

\begin{Teo}\label{teo:isofvbn/gamma}
Let $n\geq 2$. 
The groups $FWB_n$ and $FVB_n/\Gamma_2(FVP_n)$ are isomorphic to $VT_n/\Gamma_2(PVT_n)$.
\end{Teo}

\begin{proof}
We must show that $\sigma_{i+1}\sigma_i\rho_{i+1}\sigma_i^{-1}\sigma_{i+1}^{-1}\rho_i^{-1}$ and $\sigma_i\sigma_{i+1}\sigma_i\sigma_{i+1}^{-1}\sigma_i^{-1}\sigma_{i+1}^{-1}$ are elements of $\Gamma_2(PVT_n)$ for all $i=1,2,\dots,n-2$ (see Figure~\ref{quocientes}). 
From the proof of Theorem~\ref{teo:isowbn/gamma}, we can conclude that $\sigma_{i+1}\sigma_i\rho_{i+1}\sigma_i^{-1}\sigma_{i+1}^{-1}\rho_i^{-1}\in\Gamma_2(PVT_n)$ for all $i=1,2,\dots,n-2$. 
We will show that $\sigma_i\sigma_{i+1}\sigma_i\sigma_{i+1}^{-1}\sigma_i^{-1}\sigma_{i+1}^{-1}\in\Gamma_2(PVT_n)$ for all $i=1,2,\dots,n-2$, in order to do that we will use the relations given in Theorem~\ref{apvtn}, as follows:
\begin{align*}
\sigma_i\sigma_{i+1}\sigma_i\sigma_{i+1}^{-1}\sigma_i^{-1}\sigma_{i+1}^{-1}&=\rho_i\rho_i\sigma_i\sigma_{i+1}\sigma_i\sigma_{i+1}^{-1}\sigma_i^{-1}\sigma_{i+1}^{-1}\\
&=\rho_i\lambda_{i+1,i}^{-1}\sigma_{i+1}\sigma_i\sigma_{i+1}^{-1}\sigma_i^{-1}\sigma_{i+1}^{-1}\\
&=\rho_i\lambda_{i+1,i}\sigma_{i+1}\rho_{i+1}\rho_{i+1}\sigma_{i}\sigma_{i+1}^{-1}\sigma_{i}^{-1}\sigma_{i+1}^{-1}\\
&=\rho_i\lambda_{i+1,i}^{-1}\lambda_{i+1,i+2}^{-1}\rho_{i+1}\sigma_i\sigma_{i+1}^{-1}\sigma_{i}^{-1}\sigma_{i+1}^{-1}\\
&=\rho_i\lambda_{i+1,i}^{-1}\lambda_{i+1,i+2}^{-1}\rho_{i+1}\sigma_i\rho_i\rho_i\sigma_{i+1}^{-1}\sigma_i^{-1}\sigma_{i+1}^{-1}\\
&=\rho_i\lambda_{i+1,i}^{-1}\lambda_{i+1,i+2}^{-1}\rho_{i+1}\lambda_{i,i+1}^{-1}\rho_i\sigma_{i+1}^{-1}\sigma_i^{-1}\sigma_{i+1}^{-1}\\
&=\rho_i\lambda_{i+1,i}^{-1}\lambda_{i+1,i+2}^{-1}\rho_{i+1}\lambda_{i,i+1}^{-1}\rho_i\sigma_{i+1}^{-1}\rho_{i+1}\rho_{i+1}\sigma_i^{-1}\sigma_{i+1}^{-1}\\
&=\rho_i\lambda_{i+1,i}^{-1}\lambda_{i+1,i+2}^{-1}\rho_{i+1}\lambda_{i,i+1}^{-1}\rho_i\lambda_{i+2,i+1}\rho_{i+1}\sigma_i^{-1}\sigma_{i+1}^{-1}\\
&=\rho_i\lambda_{i+1,i}^{-1}\lambda_{i+1,i+2}^{-1}\rho_{i+1}\lambda_{i,i+1}^{-1}\rho_i\lambda_{i+2,i+1}\rho_{i+1}\sigma_i^{-1}\rho_i\rho_i\sigma_{i+1}^{-1}\\
&=\rho_i\lambda_{i+1,i}^{-1}\lambda_{i+1,i+2}^{-1}\rho_{i+1}\lambda_{i,i+1}^{-1}\rho_i\lambda_{i+2,i+1}\rho_{i+1}\lambda_{i+1,i}\rho_i\sigma_{i+1}^{-1}\rho_{i+1}\rho_{i+1}\\
&=\rho_i\lambda_{i+1,i}^{-1}\lambda_{i+1,i+2}^{-1}\rho_{i+1}\lambda_{i,i+1}^{-1}\rho_i\lambda_{i+2,i+1}\rho_{i+1}\lambda_{i+1,i}\rho_i\lambda_{i+2,i+1}\rho_{i+1}\\
&=\lambda_{i,i+1}^{-1}\lambda_{i,i+2}^{-1}\rho_i\rho_{i+1}\rho_{i}\lambda_{i+1,i}^{-1}\lambda_{i+2,i+1}\lambda_{i+2,i}\rho_{i+1}\rho_i\rho_{i+1}\lambda_{i+1,i+2}\\
&=\lambda_{i,i+1}^{-1}\lambda_{i,i+2}^{-1}\lambda_{i+1,i+2}^{-1}\lambda_{i,i+1}\lambda_{i,i+2}\lambda_{i+1,i+2}.
\end{align*}

So, $\sigma_i\sigma_{i+1}\sigma_i\sigma_{i+1}^{-1}\sigma_i^{-1}\sigma_{i+1}^{-1}\in\Gamma_2(PVT_n)$ for all $i=1,2,\dots,n-2$.
\end{proof}

\begin{Obs}

The \textit{virtual Gauss braid group} $GVB_n$ is the group $FVB_n$ adding the relations $\sigma_i\rho_i=\rho_i\sigma_i$, for $i=1,\dots,n-1$. Let $\pi_p\colon GVB_n\longrightarrow S_n$, where ${\pi_p}(\sigma_i)={\pi_p}(\rho_i)=\tau_i$ for $i=1,\dots,n-1$ and $\tau_i=(i,i+1)$. 
The \textit{virtual Gauss pure braid group} $GVP_n$ is defined to be the kernel of $\pi_p$. 
A presentation for $GVP_n$ was given in Proposition~5.6 of~\cite{BBD2015}, and using it we conclude that $GVB_n/\Gamma_2(GVP_n)$ is not a crystallographic group because the abelianization of the Gauss virtual pure braid group has finite order elements. 
However, concerning structural aspects, using the same techniques we can state and prove equivalent versions of Theorems~\ref{teo2},~\ref{teo3},~\ref{teo4} for the quotient group $GVB_n/\Gamma_2(GVP_n)$.
\end{Obs}

\subsection{The pure extended loop braid group}\label{subsec:ext}


We start this subsection recalling some necessary definitions to understand the objects treated here. 
The \textit{mapping class group of a $3$-manifold $M$ with respect to a submanifold $N$}, denoted by MCG$(M,\, N)$, is the group of isotopy classes of self-homeomorphisms of Homeo$(M;\, N)$, where the multiplication is determined by composition. 
The \textit{pure mapping class group of a $3$-manifold $M$ with respect to a submanifold $N$}, denoted by PMCG$(M,\, N)$, is the subgroup of elements of MCG$(M;\, N)$ that send each connected component of $N$ to itself. 
We denote by Homeo$(M;\, N^{\ast})$ the group of self-homeomorphisms of $(M,\, N)$, removing the condition of preserving orientation on $N$, with multiplication given by the usual composition.
The \textit{extended mapping class group of a $3$-manifold $M$ with respect to a submanifold $N$}, denoted by MCG$(M,\, N^{\ast})$, is the group of isotopy classes of self-homeomorphisms of Homeo$(M;\, N^{\ast})$, with multiplication determined by composition. 
The \textit{pure extended mapping class group of a $3$-manifold $M$ with respect to a submanifold $N$}, denoted by PMCG$(M,\, N^{\ast})$, is the subgroup of MCG$(M;\, N^{\ast})$ that send each connected component of $N$ to itself.

Let $n\geq 1$. Let $B^3$ denote the 3-ball and let $C=C_1 \sqcup \cdots \sqcup C_n$ be a collection of $n$ disjoint, unknotted, oriented circles, that form a trivial link of $n$ components in $\R^3$. 
The \textit{loop braid group on $n$ components}, denoted by $LB_n$, is the mapping class group MCG$(B^3,\, C)$. 
The \textit{pure loop braid group on $n$ components}, denoted by $PLB_n$, is the pure mapping class group PMCG$(B^3,\, C)$. 
Similarly, the \textit{extended loop braid group}, denoted by $LB_n^{ext}$, is the extended mapping class group MCG$(B^3,\, C^{\ast})$. 
The \textit{pure extended loop braid group}, denoted by $PLB_n^{ext}$, is the pure extended mapping class group PMCG$(B^3,\, C^{\ast})$. 
For more details about (extended) loop braid groups we refer the reader to \cite{D2017}, in particular the definitions mentioned above are available in its Section~2. 

Before stating the results about quotients of the extended loop braid groups related to crystallographic groups, we mention another interesting family of groups.
Let $n\geq 1$. The \textit{ring group $R_n$} is the fundamental group of the space of configurations of $n$ Euclidean, unordered, disjoint, unlinked circles in $B^3$. 
The \textit{untwisted ring group $U R_n$} is the fundamental group of the space of configurations of $n$ Euclidean, unordered, disjoint, unlinked circles in $B^3$ lying on planes parallel to a fixed one. 
Similarly, the \textit{pure ring group $P R_n$} is the fundamental group of the space of configurations of $n$ Euclidean, ordered, disjoint, unlinked circles in $B^3$. 
The \textit{pure untwisted ring group $P U R_n$} is the fundamental group of the space of configurations of $n$ Euclidean, ordered, disjoint, unlinked circles in $B^3$ lying on planes parallel to a fixed one. 
See \cite{BH2013} for more details about these groups.

Let $n\geq 1$. 
It was proved in \cite[Proposition~3.10]{D2017} that there is an isomorphism $PLB_n^{ext}\cong P R_n$. 
Also, we note that the pure loop braid group $PLB_n$ is isomorphic to the pure untwisted ring group $P U R_n$, and between their respective unordered versions, see \cite[Proposition~3.12]{D2017}.
Furthermore, the loop braid group $LB_n$ is isomorphic to the  welded group $WB_n$, see \cite[Corollary~6.10]{D2017}, and so are isomorphic to their respective pure versions. 

In the proof of \cite[Proposition~2.2]{BH2013} the authors proved that there is a split short exact sequence $1\to P U R_n \to P R_n \to \Z_2^n \to 1$. It follows from the isomorphisms described above that pure extended loop braid group $PLB_n^{ext}$ is isomorphic to the semi-direct product $WP_n\rtimes(\Z_2)^n$, see also the proofs of \cite[Propositions~3.12 and 4.8]{D2017}.  
In addition, the following group presentation is given for it.


\begin{Teo}[Proposition~4.8, \cite{D2017}]\label{apwpnext}
The {pure extended loop braid group} $PLB_n^{ext}$ admits the following group presentation
\begin{itemize}
\item Generators:
\subitem $\alpha_{i,j}$ for $1\le i\neq j\le n$. 
\subitem $\tau_k$ for $k=1,\dots,n$. 
\item Relations:
    \subitem $\alpha_{i,j}\alpha_{k,l}=\alpha_{k,l}\alpha_{i,j}$.
    \subitem $\alpha_{i,j}\alpha_{k,j}=\alpha_{k,j}\alpha_{i,j}$.
    \subitem $(\alpha_{i,j}\alpha_{k,j})\alpha_{i,k}=\alpha_{i,k}(\alpha_{i,j}\alpha_{k,j})$.
    \subitem $\tau_i^2=1.$
    \subitem $\tau_i\alpha_{i,j}=\alpha_{i,j}\tau_i$.
    \subitem $\tau_i\alpha_{j,k}=\alpha_{j,k}\tau_i$.
    \subitem $\tau_i\alpha_{j,i}\tau_i=\alpha_{j,i}^{-1}$.
\end{itemize}
where different letters stand for different indices. 
\end{Teo}

We have the following the short exact sequence
\begin{equation*} 
	\xymatrix{
		1 \ar[r] & WP_n \ar[r]^-{} &  PLB_n^{ext} \ar[r]^-{} &  (\Z_2)^n \ar[r] & 1,}
\end{equation*}
that provides us the following short exact sequence
\begin{equation}
	\xymatrix{
		1 \ar[r] & WP_n/\Gamma_2(WP_n) \ar[r]^-{} &  PLB_n^{ext}/\Gamma_2(WP_n) \ar[r]^-{} &   (\Z_2)^n \ar[r] & 1.}
\end{equation}
The action of $(\Z_2)^n$ on $WP_n$ is described in Theorem~\ref{apwpnext}, we rewrite it as follows
\begin{equation}\label{actionwpnext}
\tau_k\alpha_{i,j}\tau_k=
\begin{cases}
\alpha_{i,j}^{-1},\mbox{ if }k=j,\\
\alpha_{i,j},\mbox{ otherwise },
\end{cases}
\end{equation}
for all $1\le i\neq j\le n$ and $k=1,\dots,n$.

Below, we will relate a quotient of the extended loop braid group to crystallographic groups.

\begin{Teo}\label{teo1ext}
The group $PLB_n^{ext}/\Gamma_2(WP_n)$ is a crystallographic group of dimension $n(n-1)$ and whose holonomy group is $(\Z_2)^n$.
\end{Teo}

\begin{proof}
Here we use the presentation of $WP_n$ given in \cite[Corollary~3.19]{D2017}.
Then it is clear that $WP_n/\Gamma_2(WP_n)$ is isomorphic to the free Abelian group $\Z^{n(n-1)}$.
The action of $(\Z_2)^n$ on $WP_n$, given in~(\ref{actionwpnext}), induces an injective action of $(\Z_2)^n$ on $WP_n/\Gamma_2(WP_n)$. 
Thereby, the proof follows in an analogous way to the one given for Theorems~\ref{teo1} and~\ref{teo7}.
\end{proof}

Also, as given for some quotients of the virtual braid group in the last subsection, it is not difficult to prove for the group $PLB_n^{ext}/\Gamma_2(WP_n)$ equivalent versions of many results obtained for the quotient $VB_n/\Gamma_2(VP_n)$ in Section~\ref{groupvbn}.
For instance, we state below the respective versions of the results about finite order elements and its conjugacy classes. 
We note that the only possible torsion in $PLB_n^{ext}/\Gamma_2(WP_n)$ is 2.

\begin{Teo}\label{teo2ext}
\begin{enumerate}
    \item The element $\prod_{1\le r\neq k\le n}\alpha_{r,k}^{a_{r,k}}\tau_k$ have order $2$ in $PLB_n^{ext}/\Gamma_2(WP_n)$, where $a_{r,k}\in\Z$ for all $1\le r\neq k\le n$ and $1\le k\le n$.

    \item The elements $\prod_{1\le i\le j\le n}\alpha_{i,j}^{a_{i,j}}\tau_k$ and $\prod_{1\le i\le j\le n}\alpha_{i,j}^{b_{i,j}}\tau_k$ are conjugate in $PLB_n^{ext}/\Gamma_2(WP_n)$ if and only if
$$
b_{i,j}=
\begin{cases}
 a_{i,j}+2q,\mbox{ if }j=k,\\
 a_{i,j},\mbox{ if } j\neq k,
\end{cases}
$$
where $a_{i,j},b_{i,j}\in\Z$ for all $1\le i\neq j\le n$, $k=1,\dots,n$ and $q\in\Z$.

\end{enumerate}
\end{Teo}

\begin{proof}
Follows using the same ideas of the proof of Theorems~\ref{teo2} and~\ref{teo3}.
\end{proof}

\bibliography{refs}
\bibliographystyle{acm}

\end{document}